\documentclass{amsart}

\usepackage{amssymb,latexsym, amsmath, amsxtra}
\usepackage[dvips]{graphics}
\usepackage{graphicx}
\usepackage{hyperref}
\usepackage{placeins}

\theoremstyle{definition}

\theoremstyle{remark}

\numberwithin{equation}{section}
\numberwithin{figure}{section}
\numberwithin{table}{section}
\newcommand{\C}{\mathbb C}

\newcommand{\R}{\mathbb R}
\newcommand{\A}{\mathbb A}
\newcommand{\Z}{\mathbb Z}
\newcommand{\Q}{\mathbb Q}

\newcommand{\cL}{{\mathcal{L}}}
\newcommand{\mpar}[1]{}
\newcommand{\cmmt}[1]{}

\renewcommand{\Re}{\mathrm{Re}}

\DeclareMathOperator{\GL}{GL}
\DeclareMathOperator{\PGL}{PGL}
\DeclareMathOperator{\SL}{SL}

\DeclareMathOperator{\SU}{SU}
\DeclareMathOperator{\Sp}{Sp}
\DeclareMathOperator{\SO}{SO}
\DeclareMathOperator{\Orth}{O}

\newcommand{\term}[1]{\emph{\textbf{#1}}}
\newcommand{\abs}[1]{| #1 |}

\newcommand{\GammaR}{\Gamma_{\mathbb R}}
\newcommand{\GammaC}{\Gamma_{\mathbb C}}
\newcommand{\Fr}{\mbox{Fr}}
\newcommand{\Gal}{\mbox{Gal}}

\begin{document}
\title[The Landscape of L-functions: degree 3 and conductor~1]{The landscape of L-functions:\\ degree 3 and conductor~1}
\author[Farmer]{David~W.~Farmer}
\author[Koutsoliotas]{Sally Koutsoliotas}
\author[Lemurell]{Stefan Lemurell}
\author[Roberts]{David P.~Roberts}



\begin{abstract}
We extend previous lists by numerically computing approximations to many L-functions 
of degree $d=3$, conductor $N=1$, and small spectral parameters.  We sketch how previous arguments
extend to say that for very small spectral parameters there are no such L-functions.  Using the case
$(d,N) = (3,1)$ as a guide, we explain how the set of all L-functions with any fixed invariants $(d,N)$ can be viewed as a landscape of 
points in a $(d-1)$-dimensional Euclidean space.   We use Plancherel measure
to identify the expected density of points for large spectral parameters for general $(d,N)$.   
The points from our data are close to the origin and we find that they have smaller density.
\end{abstract}


\thanks{The authors thank Ralf Schmidt and Akshay Venkatesh for helpful conversations.
A portion of this work was supported by the National Science Foundation.}

\maketitle

\section{Introduction}
\label{sec:intro}

\subsection{Overview via pictures}
Consider the set $\cL_{d,N}$ of all automorphic L-functions of degree $d$ and conductor $N$,
as defined in  Section~\ref{sec:Lreview}.  
It decomposes into 
well-defined subsets of algebraic and transcendental L-functions, 
\begin{equation}
\label{decomp} \cL_{d,N} = \cL_{d,N}^{\rm alg} \coprod \cL_{d,N}^{\rm trans}.
\end{equation} The algebraic part is of particular interest because it 
conjecturally agrees with the set of  
analytically-normalized motivic L-functions of the given degree and 
conductor.  But the transcendental part is of great interest too, 
in part because for $d \geq 3$ it is much larger than the algebraic
part.  
   Tabulation efforts, such as the LMFDB, 
have focused mainly on $\cL^{\rm alg}_{d,N}$ and, via classical   
Maass wave forms, $\cL_{2,N}$.   In this paper, a sequel to  
\cite{FKL}, we aim to bring more computational attention to 
the very large sets $\cL^{\rm trans}_{d,N}$ for $d \geq 3$. 
  We focus on the first case $(d,N) = (3,1)$.  

We think about the sets $\cL_{d,N}$ visually, with L-functions being represented by 
points in regions of coordinatized $(d-1)$-dimensional Euclidean spaces. 
Each of these \term{L-points} has a natural attached multiplicity, usually $1$.  A region together with its collection of L-points 
is a \term{landscape}.  We have two organizational schemes.  One breaks $\cL_{d,N}$ into
$\lfloor(d+2)^2/4 \rfloor$ disjoint parts and places each part in its own \term{parameter landscape}.   
The other places all the points in a single \term{coefficient landscape}, where the region
is all of $\R^{d-1}$.   The names come from the nature of the coordinates, which are
real and imaginary parts of spectral parameters in the first organizational scheme and coefficients~$c_2$, \dots, $c_d$
of a polynomial in the second.

 For our case $(d,N)=(3,1)$, Figure~\ref{fig:r0r0r0} and Figure~\ref{fig:r0r1r1} each draw
a parameter landscape, while Figure~\ref{fig:Cr1plot} superimposes two similar parameter landscapes. 
The two remaining parameter landscapes for $(d,N)=(3,1)$ are theoretically obstructed
from having L-points, a special feature of the small conductor $N=1$, and hence not drawn.  Figure~\ref{eucliddata} draws
the coefficient landscape, with the contribution of each parameter landscape being indicated
via color.   It becomes computationally more difficult to find L-functions as one goes further away
from the origin in a given parameter space, and all our figures are of course partial in that they only 
show the known L-points.  Also our methods
are numerical, so some L-points could conceivably be incorrect, although
we think this is unlikely.  We encourage the reader to turn to these four figures,
because they give an accurate first sense of the content of this paper.

We are interested in understanding the landscapes both computationally, as just indicated,
and theoretically.  The foundational theoretical concepts are \term{Plancherel measure} $\mu_d$ and its 
\term{canonical approximation} $\mu_d'$.    Attention to these concepts is a significant
advance in this paper beyond the results of \cite{FKL}.  The very similar measures $\mu_d$ and
$\mu_d'$ are indicated by contour plots on our three figures of parameter landscapes.  They are not indicated directly
on the coefficient landscape because they do not need to be: $\mu'_d = P_d \, dc_2 \cdots dc_d$ 
is exactly Euclidean.  The meaning of $\mu_d$ is that it conjecturally governs the distribution of 
L-points in a precise asymptotic sense.    We combine three formulas from the 
literature to identify the \term{Plancherel constants} $P_d$.  The relevant one for our 
case is 
\begin{equation}
\label{eqn:P3}
P_3 = \frac{\sqrt{3} \zeta(3)}{128 \pi^3} \approx 0.0005246.  
\end{equation}
While our new data is only in the $(d,N) = (3,1)$ case, our text works with general $(d,N)$ as it   
is aimed at supporting future numerical investigations.

\subsection{Content of the sections}
Section~\ref{sec:Lreview} defines the sets $\cL_{d,N}$, reviews standard material about L-functions,
and defines the $\lfloor (d+2)^2/4 \rfloor$ receiver regions for the parameter landscapes.  
Section~\ref{sec:Lfunctions} populates the landscapes with L-points, following the methods of \cite{FKL}.  
Section~\ref{sec:free} offers reminders that each individual L-function is interesting and sketches
how one can prove that parameter landscapes where the conductor $N$ is small enough 
 have an L-point-free region near their origin.   

Section~\ref{sec:coefficient} discusses the purely algebraic passage from parameter landscapes to the coefficient
landscape by means of symmetric functions.   Section~\ref{sec:plancherel} gives the Plancherel theory including
the identification of the constants $P_d$.  Our data gives points at roughly half the asymptotic density 
$P_d$.   We explain by comparison with the well-understood $d=2$ case
that this is an expected shortfall; there should be as yet undiscovered secondary terms
which will make theory better match the data.  Finally Section~\ref{sec:complements} gives two independent 
complements.  One compares our data with the $p$-adic Plancherel measures on the class space
of $\SU_3$, finding again the need for secondary terms.   The other discusses the
extent to which transcendental L-functions outnumber algebraic L-functions
in degrees $\geq 3$.

\subsection{Availability of data}  The L-functions in $\cL_{3,1}$ have the form 
\begin{equation}
\label{eqn:eulerproduct3}
L(s) = \prod_{{\rm primes \; } p} \frac{1}{1-a_p p^{-s} + \overline{a}_p p^{-2s} - p^{-3s}}.
\end{equation}
For each drawn L-point, its coordinates in the landscape and the complex numbers
$a_p$ to their computed precision are available on GitHub at~\cite{FKLrep}.
Earlier data, dating back to \cite{FKL, Bi, Bo}, is already conveniently
available on the LMFDB for $(d,N) \in \{(3,1),(3,4),(3,9),(4,1)\}$. 
We plan to systematically extend this section of the LMFDB,
because much more is within computational reach.

%

\section{Background on L-functions}
\label{sec:Lreview} 
         Here we review L-functions with an initial focus on the fundamental
 decomposition of all automorphic L-functions into the two types,
 \begin{equation}
 \cL= \cL^{\rm alg} \coprod \cL^{\rm trans}.
 \end{equation}
 Whether an automorphic L-function is algebraic or transcendental depends only
 on the $\Gamma$-factor in its functional equation and our second main
 focus is the space $X_d$ of possible $\Gamma$-factors in degree $d$,
 and various related spaces.    
         

      

\subsection{The set of automorphic L-functions $\cL$}
As a terminological catch-all, define an \term{L-function} to be a Dirichlet series with 
Euler product converging in some right half plane:
\begin{equation}
  \label{L}
  L(s) = \prod_{{\rm primes \; } p} \frac{1}{f_p(p^{-s})} = \sum_{n=1}^{\infty} \frac{a_n}{n^s};
 \end{equation}
 here the $f_p(x) = 1 + \cdots$ are required to be polynomials in $\C[x]$ having a maximal 
 \term{degree} $d$.  For example, the Riemann zeta function $\zeta(s)$ has $f_p(x)=1-x$ 
 for all $p$ and degree $1$.  The class of 
 L-functions just defined is not at all interesting in its entirety.  However
 it allows us to conveniently talk about automorphic, motivic, and
 our experimental L-functions.  
 
 An \term{automorphic L-function} for this paper is an L-function associated
 in the standard way to a balanced tempered-at-$\infty$ cuspidal automorphic
 representation of the adelic group $GL_d(\A)$.   We do not 
 need to enter into the very complicated automorphic theory, as we will just be
 using known properties of these L-functions.  A list of known
 properties which moreover conjecturally characterizes automorphic
 L-functions is given in \cite{FPRS}.  Here we review what we need.  
 
 Associated to an automorphic L-function is a \term{conductor} $N \in \Z_{\geq 1}$. Also 
 associated is 
 a \term{central character} $\chi$, which is a Dirichlet character of conductor $N$.
 So $\chi: \Z/N \rightarrow \C$ is a multiplicative homomorphism, making 
 $\chi(p) = 0$ for the primes $p$ dividing $N$.   The L-functions coming
 from $GL_d(\A)$ have degree $d$, as
 \[
 f_p(x) = 1 - a_p x + \cdots + (-1)^d \chi(p) x^d.
 \]
 The Riemann zeta function is automorphic, and has a meromorphic
 continuation to the whole plane with a unique pole at $s=1$.  All
 other automorphic L-functions are entire.  The completed L-function and its functional 
equation take the form
\begin{align}
\label{Ldef} \Lambda(s) & :=\mathstrut L(s) \, N^{s/2} 
\prod_{j=1}^{d_1} \GammaR(s+ \mu_j)
\prod_{k=1}^{d_2} \GammaC(s+ \nu_k) \\
\label{funct1}  & = \varepsilon \overline{\Lambda}(1-s).
\end{align}
Here the bar indicates Schwarz reflection, so that $\overline{L}(s) = \overline{L(\overline{s})} = \sum_n \overline{a}_n n^{-s}$ is the \term{dual} of $L$.  
Also $\GammaR$ and $\GammaC$ are the normalized $\Gamma$-functions,
\begin{equation}
\GammaR(s) = \pi^{-s/2} \Gamma(s/2),
\ \ \ \ \ \ \
\ \ \ \ \ \ \
\ \ \ \ \ \ \
\GammaC(s) = (2\pi)^{-s} \Gamma(s).
\end{equation}
The double product of $\Gamma$-functions in \eqref{Ldef} is the \term{$\Gamma$-factor} of 
$L$, to be discussed in more detail in the next subsection.   The \term{sign} $\varepsilon$ is 
a complex number on the unit circle. The $d_i \in \Z_{\geq 0}$ satisfy $d_1+2d_2=d$
and form the {\bf signature} $(d_1,d_2)$ of $L$.

\subsection{$\Gamma$-factors}
\label{ssec:terminology}
Much of this paper is driven by the nature of $\Gamma$-factors, and so the terms introduced in 
this section are particularly important.     The \term{spectral parameters} $\mu_j$ and $\nu_k$ in \eqref{Ldef}
are written in terms of their real and imaginary parts as
\begin{align*}
\label{realimag} \mu_j & = \delta_j + i \lambda_j, & 
\nu_k &=  \kappa_k + i \beta_k.
\end{align*}
They are constrained by the \term{balanced} condition in our definition of 
automorphic L-function, $\sum \lambda_j + 2 \sum \beta_k = 0$.  
They are constrained also by the \term{tempered-at-$\infty$} condition,
which requires $\delta_j \in \{0,1\}$ and $\kappa_k \in \{\frac12, 1, \frac32, 2, \dots\}$.
There is no ambiguity in $d_1$ and $d_2$ coming from the duplication formula $\Gamma_\C(s) = \Gamma_\R(s) \Gamma_\R(s+1)$ because 
$\kappa_k=0$ is disallowed.   
 Important for us is the \term{refined signature} $(d_+,d_-,d_2)$,
 where $d_+$ and $d_-$ are the number of $j$ with $\delta_j=0$ and $\delta_j=1$ respectively.
Henceforth, we write $\lambda_{d_1+k}$ rather than $\beta_k$.  
To remove ambiguities, we require
that the $\mu_j$ and $\nu_k$ are both increasing with respect
to the standard lexicographical order on $\C$. 

The \term{$\Gamma$-type} of a $\Gamma$-factor is the list of $\delta_j$ followed
by the list of $2\kappa_k$.  For example,
a degree seven L-function might have $\Gamma$-type $(0,0,1,1,1;5)$.
Following LMFDB conventions, we normally write 
$\Gamma$-types by prefixing each $\delta_i$ with an $r$ and each
$2 \kappa_k$ with a $c$, as in $r0r0r1r1r1c5$.   A $\Gamma$-type is 
called \term{even} if all $2 \kappa_k$ are even.  It is called
\term{odd} if $d_1=0$ and all $2 \kappa_k$ are odd.  Thus 
our example $\Gamma$-type is neither odd nor even.

There is an important constraint on the spectral parameters
in terms of the central character:
\begin{equation}
\label{eqn:chi}
\chi(-1) = (-1)^{\sum \delta_j +  \sum(2\kappa_k+1)}.
\end{equation}
Always $\chi(-1) = \epsilon_\infty^2$, where $\epsilon_\infty \in \{1,i,-1,-i\}$ depends
multiplicatively on the individual Gamma factors in $\eqref{Ldef}$, in the way 
presented in Table~\ref{onetwothree}.  
%
%
 There are two simplifications when $N=1$.  First, $\chi(-1)=1$, allowing one to
 disregard half the possible spectral parameters, namely those that fail \eqref{eqn:chi}.
 Second $\varepsilon = \epsilon_\infty$; this replaces an unknown quantity with a known
 one in the searches described in \S\ref{ssec:searchmethod}.
 
  \subsection{The decomposition $\cL = \cL^{\rm alg} \coprod \cL^{\rm trans}$} 
  By definition, an automorphic L-function is \term{algebraic} if its spectral parameters 
  are real and its $\Gamma$-type has a parity.   The remaining 
  automorphic L-functions are \term{transcendental}.    The next two paragraphs explain
  how this innocuous-looking decomposition is expected to
  be a very sharp dichotomy.  
  
  A \term{motivic L-function} $L$ comes from an elaborate process 
  that starts by counting points on a fixed algebraic variety over varying finite 
  fields.    It has an associated $\Gamma$-factor coming from 
  Hodge theory.   We assume a primitivity condition corresponding
  to focusing attention on an irreducible motive.  After switching to the analytic normalization,
  a change of variables of the form $s \rightarrow s+\frac{w}{2}$
  for $w$ an integer, a fundamental and widely-believed conjecture says that 
  motivic L-functions are automorphic.  By the nature 
  of their $\Gamma$-factors, they would have to be in
  $\cL^{\rm alg}$.   A~conjectural converse says 
  that all L-functions in $\cL^{\rm alg}$ come from
  motives in this way. 
  
  For a motivic L-function, even in its analytic normalization, all the coefficients 
  $a_n$ are algebraic numbers.   In contrast, for a transcendental 
  L-function one expects that all the nonzero $\lambda_j$ and almost all the nonzero $a_p$ 
  are transcendental numbers.   Assuming that all algebraic L-functions are motivic,
  this conjecture translates to a contrast between 
  $\cL^{\rm alg}$ and $\cL^{\rm trans}$.  The contrast also 
  explains the names of the summands. 
    
  We refer the reader again to \cite{FPRS} for more details on everything
  discussed so far.  One subtlety worth highlighting concerns
  Axioms 4a and 4b on the list of axioms in Section 2.1 of \cite{FPRS}.  We 
  are imposing Axiom 4a by our ``temperedness-at-$\infty$" 
  condition, because it is necessary for our points-in-landscapes
  presentation.  The Selberg conjecture says that
  we would not enlarge our set $\cL$ if we removed
  this condition.  We are not imposing Axiom 4b, which 
  corresponds to a ``temperedness at all finite primes $p$'' condition, 
  because it is not necessary for our presentation. 
  Here the Ramanujan conjecture says that we would not
  decrease our set $\cL$ if we added this condition.  
  
  \subsection{The constructional problem}
  Consider now the problem of explicitly writing down
  algebraic L-functions, say in an approximate form
  $\sum_{n=1}^c a_n n^{-s}$.   One can do this easily
  by purely automorphic methods only in very
  limited circumstances.  For example, consider
  automorphic L-functions with $\Gamma$-factor 
  $\Gamma_\C(s+1/2)^g$ and
  $\chi$ the trivial character.  
   Then for $g=1$, automorphic methods are successful: one
  can find a corresponding holomorphic newform
  of weight $2$ on some $\Gamma_0(N)$.   
  For $g \geq 2$, it is not immediately
  clear by automorphic methods that $\Gamma_\C(s+1/2)^g$
  even actually arises, and explicitly computing a
  particular $\sum_{n=1}^c a_n n^{-s}$ even for small cutoff $c$ seems problematic.
  In contrast, writing down a genus $g$ curve 
  over $\Q$ and computing many $a_n$ for it is much
  more straightforward.  
  
  Consider now the same constructional problem but
  now for transcendental  L-functions.  
  Purely automorphic methods seem even less
  promising in this context.   In this paper, we do
  what seems like the best that can be currently
  done: we produce \term{experimental L-functions}
  $\prod_{p} f_p(x)^{-s}$ having only finitely many 
  non-trivial factors.   Their purpose is to approximate actual 
  automorphic L-functions.  At present,
  like in typical instances of using motivic
  L-functions as in the last paragraph, 
  the evidence for the existence of a corresponding
  automorphic L-function is strong, but not conclusive. 
  
  Having carefully made distinctions, we will now generally
  be more brief with our language:  ``L-function'' in theoretical contexts will 
  mean automorphic L-function;  ``L-function" in the context
  of our searches will mean experimental L-function.  
  One would like to populate our landscapes with actual
  automorphic L-points, but we are of course populating
  them with experimental L-points.

\subsection{Spaces of generalized $\Gamma$-factors}  
\label{ssec:spaces} 
The material of this subsection is completely elementary.  However it is not standard, but 
rather is particular to our concerns in this paper.  To define the Euclidean 
regions receiving L-points, it is necessary to allow 
%
%
 \term{generalized $\Gamma$-factors} by relaxing the half-integrality condition 
on the $\kappa_k$ to $\kappa_k \in (0,\infty)$.  
Let  $X'_{d_-,d_+,d_2}$ be the set
of $\Gamma$-factors with refined signature $(d_-,d_+,d_2)$.  
Because of our ordering convention, the $\lambda_j$ and $\kappa_k$ are well-defined 
real-valued functions on $X'_{d_-,d_+,d_2}$.  
Because of the balanced condition, the functions 
$(\lambda_1,\dots,\lambda_{d_1+d_2-1},\kappa_1,\dots,\kappa_{d_{2}})$
render $X'_{d_-,d_+,d_2}$ a full-dimensional cone in 
$\R^{d-1}$ and thus a connected topological space.  
Let $X'_{d_1,d_2}$ be the disjoint union of the 
$X'_{d+,d-,d_2}$ with $d_1 = d_+ + d_-$.  
Let $X'_{d}$ likewise be the disjoint union of the 
$X'_{d_1,d_2}$ with $d = d_{1}+2d_2$.   
So clearly $X'_{d_1,d_2}$ has $d_1+1$ connected 
components. As reported in the introduction, the number of 
connected components of $X_d'$ works out to the  ``quarter
square'' $\lfloor (d+2)^2/4 \rfloor$.   
As an alternative notation for $X'_{d_+,d_-,d_2}$ 
we use repetition, as in $X'_{2,3,1} = X'_{++---c}$.
The possible spaces $X'_{d_+,d_-,d_2}$ with $d \leq 3$ 
are listed in this alternative notation in Table~\ref{onetwothree}.    

\begin{table}[htb]
\[
{\renewcommand{\arraycolsep}{2.5pt}
\begin{array}{ccc|cc| l | cc}
d& \chi(-1) & \epsilon_\infty & \gamma & \delta& \mbox{Automorphic Source} & \gamma & \dim X'_\delta  \\
\hline
1 & 1 &1&  r0 & +  & \mbox{Even Dirichlet characters} & 0& 0  \\
   & -1 &i&  r1 & - &  \mbox{Odd Dirichlet characters} & 0 & 0 \\
\hline
2 & 1 &1& r0r0 & ++ & \mbox{Even weight 0 Maass forms} & 1 & 1\\ 
   & -1   &i& r0r1 & +- & \mbox{Weight 1 forms} &1 & 1   \\
   & 1  &-1& r1r1 & --  & \mbox{Odd weight 0 Maass forms} & 1 & 1  \\
   & (-1)^k   &i^k& c2\kappa & c & \mbox{Weight $k=2\kappa+1$ forms} & 0 & 1  \\
\hline
3 & 1 &1& r0r0r0 & +++ & \mbox{Spherical Maass forms} & 2 & 2 \\
   & -1 &i& r0r0r1 & ++-   && 2 & 2 \\
   & 1 &-1& r0r1r1 & +--  && 2 & 2 \\
   & -1 &-i& r1r1r1 &  --- &&  2 & 2\\
   &(-1)^k  &i^k& r0c2\kappa & +c && 1 & 2 \\
   &(-1)^{k+1}   &i^{k+1}& r1c2\kappa & -c && 1 & 2 \\
\end{array}
}
\]
\caption{\label{onetwothree} $\Gamma$-types $\gamma$ and refined signatures
$\delta = (d_+,d_-,d_2)$ indexing inclusions $X_\gamma \subseteq X'_\delta$ of connected
spaces in degree $d \leq 3$}
\end{table}

The subspaces of all the above 
spaces consisting of actual 
$\Gamma$-factors are denoted
by removing the prime.  The
connected components of
$X_d$ are indexed by their $\Gamma$-types, as discussed 
above.
All the components of $X_{d_1,d_2}$ have 
dimension $d_1+d_2-1$, there being infinitely many 
components exactly when $d_2>0$.  One of the 
virtues of LMFDB notation for $\Gamma$-types 
is that it lets one easily distinguish the connected spaces
$X'_{d_+,d_-,d_0}$ from the connected components 
of their subspaces $X_{d_+,d_-,d_2}$.   
As an example, $r0r0r1r1r1c5 \subset X'_{++---c}$ is
an inclusion of a connected five-dimensional space into 
a connected six-dimensional space. In general, the notational transition is made by $r0 \rightarrow +$,
$r1 \rightarrow -$, and replacing all the  $2\kappa_k$'s by 
$c$'s.  All possibilities for $d \leq 3$ are listed in 
Table~\ref{onetwothree}.

\subsection{Sources} 
 \label{ssec:sources}
 While we do not need the full automorphic theory, we will reference automorphic
 forms in classical contexts.  Very simply, it is known that $\cL_{1,N}= \cL^{\rm alg}_{1,N}$ is the set 
 of Dirichlet L-functions of conductor $N$.    In the next dimension, 
$\cL_{2,N}$ is conjecturally the set $L(f,s)$ of L-functions coming
from newforms on $\Gamma_1(N)$.   Holomorphic forms 
and Maass forms with Laplacian eigenvalue $\frac{1}{4}$ give
L-functions in $\cL_{2,N}^{\rm alg}$.  Maass forms with Laplacian eigenvalue
$\lambda > \frac14$ give L-functions in $\cL_{2,N}^{\rm trans}$.   
There are not expected to be Maass forms with Laplacian
eigenvalue $<\frac14$ because their L-functions would be counterexamples to
the Selberg conjecture; if such L-functions exist, they would not 
be in our $\cL_{2,N}$.   For $d \geq 3$, 
one is already out of the classical setting.  
Table~\ref{onetwothree} summarizes the situation.  


\section{The parameter landscapes for L-functions with $(d,N)=(3,1)$}
\label{sec:Lfunctions}
Here we draw the three pictures promised in the introduction and discuss
how we found the L-functions corresponding to the L-points in the pictures.  

\subsection{The landscape of L-functions with $\Gamma$-type $r0r0r0$ in 
the cone $X_{+++}$} 
\label{ssec:r0r0r0} This is the 
landscape studied in \cite{FKL} and we roughly double the
number of L-points previously available in it.  As explained in \S\ref{ssec:spaces}, the L-functions
considered here have $\Gamma$-factor of the form 
%
%
\begin{equation}\label{eqn:r0r0r0Gamma}
\GammaR(s + i \lambda_1) \GammaR(s + i \lambda_2)  \GammaR(s + i \lambda_3),
\end{equation}
with $\lambda_j\in \R$ and $\lambda_1 + \lambda_2 + \lambda_3 = 0$.  
The general prescription of \S\ref{ssec:terminology}, useful for uniform theoretical
discussions, is in this case to require $\lambda_1 \leq \lambda_2 \leq \lambda_3$,
and use $\lambda_1$ and $\lambda_2$ as coordinates.   Instead, we will be consistent 
with \cite{FKL} and order by $\lambda_3 \leq \lambda_2 \leq \lambda_1$, 
still using $\lambda_1$ and $\lambda_2$.    This difference in
conventions is unimportant:  it is a question of choosing a fundamental domain 
for a Euclidean plane with an action by the symmetric group $S_3$.   

\begin{figure}[!htb]
\includegraphics[width=0.98\textwidth]{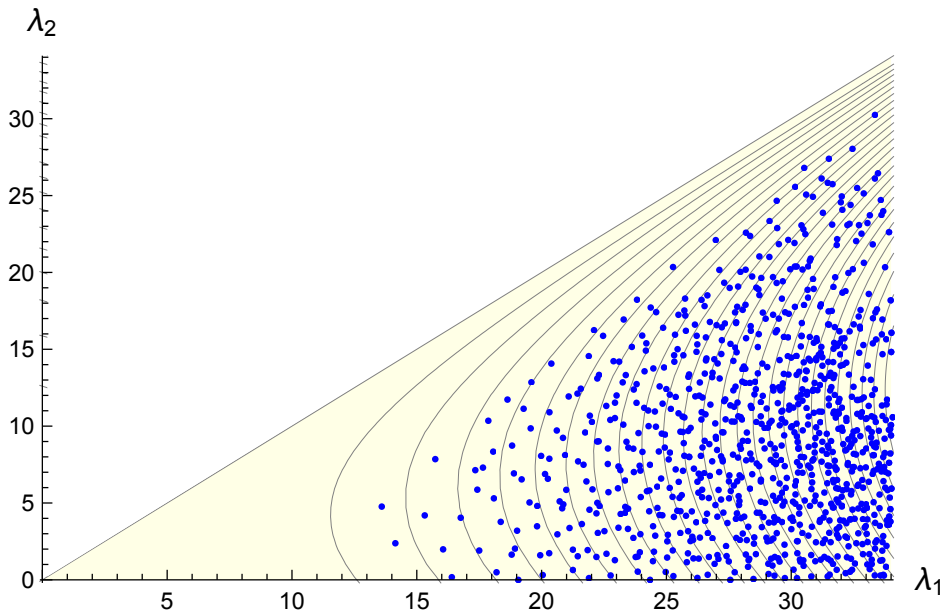}
\caption{\sf
The L-points $(\lambda_1,\lambda_2)$ for L-functions with $\Gamma$-type $r0r0r0$,
having $\Gamma$-factor
\eqref{eqn:r0r0r0Gamma} and conductor~1.
} \label{fig:r0r0r0}
\end{figure}

There is a second redundancy not emphasized in the previous section: an
L-function and its dual determine each other by complex conjugation.  To 
remove this redundancy, again following \cite{FKL}, we require 
$\lambda_2 \geq 0$.  Accordingly, Figure~\ref{fig:r0r0r0} draws only
half of the landscape; the other half is obtained by reflecting
in the $\lambda_1$-axis.    So the shaded region is a fundamental
domain for the action of a dihedral group $D_6 = S_3 \times S_2$.  
It is naturally one-twelfth of its plane, but our 
symmetry-breaking choice to emphasize $\lambda_1$ and 
$\lambda_2$ obscures this fact.




There are six L-points in Figure~\ref{fig:r0r0r0} on the $\lambda_1$-axis.  
They arise as the symmetric squares of L-functions of the 
classical Maass forms in appearing as the six rightmost colored L-points
in Figure~\ref{2pict}.  The coordinates $(\lambda_1,0)$ in Figure~\ref{fig:r0r0r0}
and $(\lambda)$ in Figure~\ref{2pict} are related by $\lambda_1 = 2 \lambda$.
%
%
 It is known that 
the non-algebraic self-dual L-functions in $\cL_{1,3}$ arise in
this way.  


A careful look at Figure~\ref{fig:r0r0r0} reveals striations running along
lines of approximate slope $0$ and $1$.  It would be interesting to 
find a theoretical explanation for this pattern.







\subsection{The landscape of L-functions with $\Gamma$-type $r0r1r1$ in 
the cone $X_{+--}$} 
\label{ssec:r0r1r1}

Our second landscape is populated with L-functions with
$\Gamma$-factor of the form
\begin{equation}\label{eqn:r0r1r1Gamma}
\GammaR(s + i \lambda_1) \GammaR(s + 1 + i \lambda_2)  \GammaR(s + 1 + i \lambda_3).
\end{equation}
Here $\lambda_2$ and $\lambda_3 = -\lambda_1-\lambda_2$ are playing the same
role and so we normalize by requiring $\lambda_2 \geq \lambda_3$.   The bottom
boundary of the drawn region comes from the equation $\lambda_2=\lambda_3$,
i.e.\ $2 \lambda_1+\lambda_2=0$.   As $\lambda_1$ is playing a different role 
from the other two $\lambda_i$, we can remove the redundancy of duality by restricting
attention to $\lambda_1 \geq 0$.    So the highlighted region
is a fundamental domain for a group of the form $C_2 \times C_2$,
and so really should be thought of as one-quarter of the full plane. 

\begin{figure}[htp]
\includegraphics[width=0.98\textwidth]{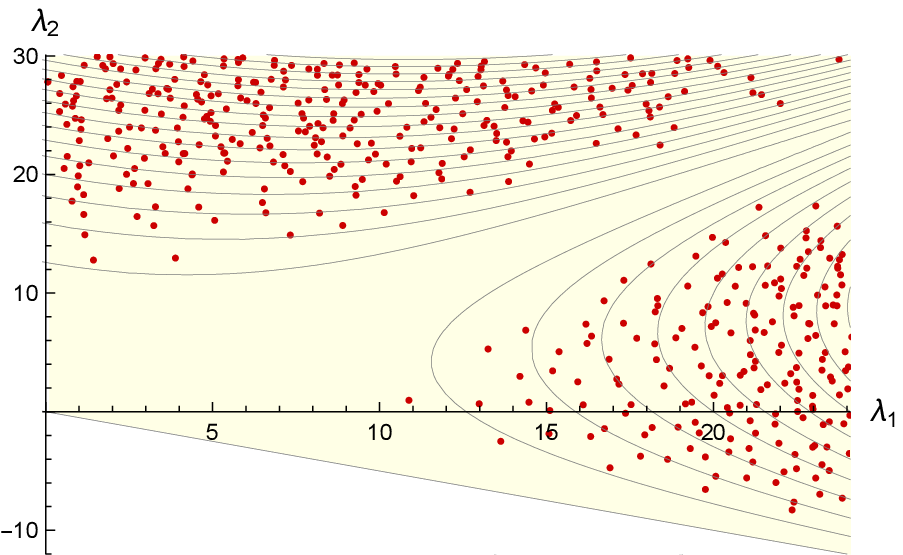}
\caption{\sf
The L-points $(\lambda_1,\lambda_2)$ for L-functions with $\Gamma$-type $r0r1r1$ having $\Gamma$-factor
\eqref{eqn:r0r1r1Gamma} and conductor~1.
} \label{fig:r0r1r1}
\end{figure}

In contrast with the previous subsection, where all L-functions 
have sign $1$, here all L-functions have sign $-1$, by Table~\ref{onetwothree}.  These
signs are important in our search for L-functions.  However, 
the sign being $-1$ does not at all force central vanishing of the L-functions, since
they are not self-dual.   Applying general principles to 
our cases here, we would expect that $L(1/2) \neq 0$ for
all any $L$ in any $\cL^{\rm trans}_{3,N}$.






\subsection{The landscape of $r0c2\kappa$ lines in $X'_{+c}$ and 
the landscape of $r1c2\kappa$ lines in $X'_{-c}$}  Our third case 
illustrates typical behavior for arbitrary $d$ better than the 
first two cases.  Here the L-functions have $\Gamma$-factors
of the form
\begin{equation}\label{eqn:rdeltakappa}
\GammaR(s + \delta + i \lambda) 
\GammaC(s + \kappa - i \lambda/2).
\end{equation}
We work in the $(\kappa,\lambda)$ plane.  
For $\tau \in \{+,-\}$, the set $X'_{\tau c}$  is identified with
the half-space where $\kappa>0$.   Following
the general prescription, a point in $X'_{+c}$
is a generalized $\Gamma$-factor \eqref{eqn:rdeltakappa}
with $\delta=0$; similarly, points in 
$X'_{-c}$ index generalized $\Gamma$-factors
\eqref{eqn:rdeltakappa} with $\delta=1$.  
To remove the redundancy given by duality,
we draw only the quadrant where also
$\lambda \geq 0$.  

\begin{figure}[htp]
\includegraphics[width=0.98\textwidth]{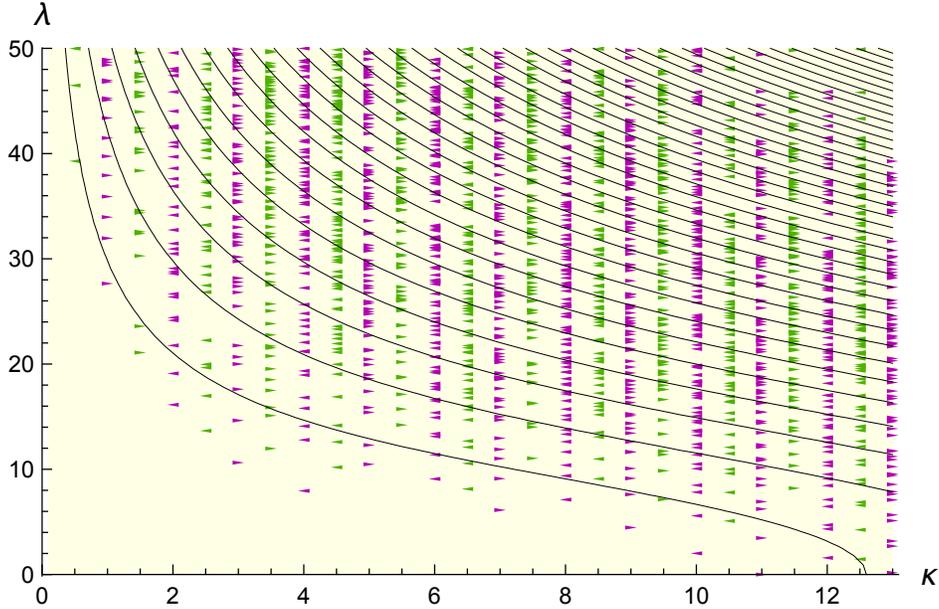}
\caption{\sf
The L-points $(\kappa,\lambda)$ for L-functions with 
$\Gamma$-type $r\delta c2\kappa$, having $\Gamma$-factor
\eqref{eqn:rdeltakappa} and conductor~1.
The green triangles correspond to $\delta=0$ and the purple triangles
have $\delta=1$.  The triangles are pointing left if the sign
of the functional equation is $-1$, and pointing right if the sign is~$+1$.
} \label{fig:Cr1plot}
\end{figure}

The spaces $X_{\tau c}$ of actual $\Gamma$-factors come from imposing 
the condition that $\kappa \in \frac{1}{2} \Z_{\geq 1}$.
Because of the central character condition \eqref{eqn:chi} and the identification
of $\chi(-1)$ on Table~\ref{onetwothree}, 
L-points in $X_{-c}$ can only have integral $\kappa$
coordinates while L-points in $X_{+c}$ can
only have half-integral $\kappa$.  It is appropriate to superimpose the 
two landscapes in one figure, because the formula \eqref{planchform} for Plancherel measure 
$\mu$ on $X_{\tau 0}$ and the formula \eqref{planchform2} for approximate
Plancherel measure $\mu'$ on $X'_{\tau 0}$ are 
independent of $\tau \in \{-,+\}$.   This agreement
of measures holds always between the spaces indexed 
by $(a,b,d_2)$ and $(b,a,d_2)$.   
Thus $X_{+++}$ from \S\ref{ssec:r0r0r0} would be naturally
superimposed with $X_{---}$ and $X_{+--}$ from
\S\ref{ssec:r0r1r1} would be naturally superimposed with 
$X_{-++}$.   To repeat a point from the introduction,
it is only because $N=1$ that each of the latter spaces
do not have L-points.  

\subsection{The algebraic L-points}
\label{ssec:algebraicL}
There is only one algebraic L-point in our three figures, 
it being $(11,0)$ in Figure~\ref{fig:Cr1plot}.  
However for the analogous figures 
with increasing $N$, algebraic L-points start
to appear as follows.  

Any space $X_{d_+,d_-,0}$ has exactly one
algebraic point, the origin $(0,\dots,0)$.  
All the L-functions mapping to $(0,\dots,0)$ 
are expected to come from the
continuous representations of 
$\mbox{Gal}(\overline{\Q}/\Q)$ into $\GL_d(\C)$ 
for which complex conjugation has eigenvalues  
$1$ and $-1$ with the respective multiplicities 
$d_+$ and $d_-$.    The LMFDB catalogs
many examples, with $d_+-d_-$ given in the column $\chi(c)$.   
Self-dual L-functions can be searched by making the Frobenius-Schur 
indicator $1$.  For the four spaces $X_{+++}$, 
$X_{++-}$, $X_{+--}$, and $X_{---}$, the smallest conductors appearing are 
respectively $1957$, $283$, $229$, and $2828$.  The images
of the representations are the groups $S_4$, $S_4$, $S_4$, and $S_4 \times S_2$,
and the L-functions all come from Dirichlet twists of symmetric squares of rank two L-functions.
Non-self-dual L-functions have Frobenius-Schur indicator $0$.   For the first three spaces, the 
smallest conductors appearing are respectively $8473$, $1228$, and $8100$, from groups
$C_3 \wr S_3$, $(C_3^2\!:\!C_3)\!:\!C_3$, and $(C_3^2\!:\!C_3)\!:\!C_3$.  As of this
writing, the LMFDB does not have an example with eigenvalues $---$, 
but many can be obtained via twisting a $+++$ example with 
an odd Dirichlet character. 




For the spaces $X_{\tau c}$, the algebraic points are $(\kappa,0)$ for 
$\kappa$ a positive integer.  In direct analogy with 
a paragraph in \S\ref{ssec:r0r0r0}, modular forms of weight $k=\kappa+1$ 
give a large supply, via symmetric squares and subsequent Dirichlet 
twists.   The aforementioned algebraic point $(11,0)$ in 
Figure~\ref{fig:Cr1plot} comes from the Ramanujan
function $\Delta \in S_{12}(1)$.   
%
%
%
%
The very low L-point on the $\kappa=13$ line
has $\lambda \approx 0.1660$.  
The fact that the corresponding L-function is not self-dual is seen
clearly also in its first coefficients:   $a_2 \approx 0.0256 - 0.7997i$ 
and $a_3 \approx -0.5418 + 1.2019 i$ are far from real.  
\S\ref{ssec:algebraicpart} discusses the part of $\cL^{\rm alg}_{3,N}$ 
with $\kappa > 0$ that does not come from modular
forms.


\subsection{The search method}
\label{ssec:searchmethod}
     Our method for locating L-points in landscapes was explained in
\cite{FKL}.  We are working with more general $\Gamma$-factors, 
but the basic method has not changed.  We give a quick summary
here.  

The starting point is the ``approximate functional equation''
    \begin{align}
         \label{eq:formula}
         \Lambda(s) g(s) =
          Q^s 
           \sum_{n=1}^{\infty} \frac{a_n}{n^s} f_1(s,n) 
         + \varepsilon Q^{1-s} 
                \sum_{n=1}^{\infty} \frac{\overline{a}_n}{n^{1-s}} f_2(1-s,n).
    \end{align}
Here $Q=\sqrt{N}$ and $f_1$ and $f_2$ are integrals involving a test function $g(s)$ and the $\Gamma$-factors. 
The main idea is that the  test function $g(s)$ is a free parameter.  Thus, by
evaluating the L-function in different ways but at the same point $s$, one obtains
an equation in the Dirichlet coefficients~$a_n$.  Besides the linear equations 
coming from \eqref{eq:formula} there are also non-linear equations coming from 
the Euler product \eqref{eqn:eulerproduct3}, reducing the set of unknowns
to the prime-indexed coefficients~$a_p$.   If the functional equation
parameters correspond to an actual L-function, then that system of
equations will be consistent.  If the parameters do not correspond to
an actual L-function, then the ability to quantify the inconsistency of
the system can be turned into a method of searching for the correct
parameters.  The solution to the consistent system will be the 
correct $a_p$'s.

To illustrate the idea, we sketch how we found the L-point with coordinates
$(\kappa,\lambda) = (2.5,16.97618\dots)$ on Figure~\ref{fig:Cr1plot}.    We work 
with an associated real-valued Z-function with  $Z(t)$
directly related to $L(\frac{1}{2}+it)$ as discussed in \S\ref{ssec:threeL} below.  
Consider a hypothetical (but actually non-existent) L-point with exact spectral parameters say 
$(\kappa,\lambda)=(\frac{5}{2},16.97)$.  Equation \eqref{eq:formula} 
with $g(s)=1$ says
\begin{align*}
Z(5)=\mathstrut & 2.44 - 1.91 a_2^r - 16.65 a_2^i 
- 1.01 a_3^r - 1.27 a_3^i - 0.091 a_4^r - 0.090 a_4^i +
\cdots\cr
& - 0.000023 a_7^r - 0.000023 a_7^i + \cdots -
5.0 \times 10^{-10} a_{11}^r - 6.2 \times 10^{-10} a_{11}^i + \cdots.
\end{align*}
With $g(s)=e^{is/2}$, the equation \eqref{eq:formula}  says
\begin{align*}
Z(5)=\mathstrut & 
-0.145 + 0.42 a_2^r - 1.54 a_2^i 
+ 0.35 a_3^r + 0.29 a_3^i - 0.042 a_4^r + 0.032  a_4^i +
\cdots\cr 
& + 0.000022  a_7^r + 0.000012 a_7^i + \cdots +
3.5 \times 10^{-10} a_{11}^r - 1.1 \times 10^{-0} a_{11}^i + \cdots .
\end{align*}
The notation $a_n^r$ and $a_n^i$ refers to the real and imaginary parts
of the Dirichlet coefficient $a_n$, respectively.  The 
coefficients come from the integrals $f_1$ and $f_2$ 
and can easily be obtained to hundreds of digits of
precision.  

%

The equality of the two expressions for $Z(5)$ gives one equation in the  
Dirichlet coefficients.  Choosing other test functions likewise gives
other equations.    The rapid decrease in the contributions of the
coefficients enables one to truncate the expression to obtain an
equation in finitely many variables. 
We can create more equations than unknowns and then solve a sub-system
with equal number of equations and unknowns.  Typically that system has
a solution, and often it has several (recall that it is a nonlinear system
because of the relations coming from the Euler product).
Plugging these solutions into the
unused equations (which we call ``detectors'')
gives a measure for the consistency of the overdetermined
system.

We repeat the process for another hypothetical (but also non-existent) L-point,
say with $(\kappa,\lambda)=(\frac{5}{2},16.98)$.  By comparing the residuals
from the detectors at the two hypothetical L-points, we can use the
secant method to determine a more accurate estimate for the
coordinates of an actual L-point.
When it is successful, our values for $\lambda$ and the $a_p$ are improved,
typically obtaining one or two additional decimal places of accuracy on
each iteration.  
We generally stop when we are confident that we have computed 
the spectral parameters to at least a dozen
digits of accuracy.   In this example, we stopped
at 
\begin{equation}
\label{Lcoeffs}
\begin{array}{rcrlrl}
\lambda & \approx & & \!\!\!16.9761877662491 \!\!\!\!\!\!\\
a_2 & \approx & -\!\!\!&0.063424433675 \!\!\! & +\!\!\!&0.10988282023 \, i \\
a_3 & \approx & &0.63801509281 & +\!\!\!&0.50484439411 \, i \\
a_5 & \approx & &0.9093205878 & -\!\!\!&0.1141574688  \, i\\
a_7 & \approx & &0.68168545 & +\!\!\!&0.02028876 \, i\\
a_{11} & \approx & - \!\! \!&0.33243 & +\!\!\!&0.12179 \, i\\
a_{13} & \approx & -\!\!\!&0.46 &  -\!\!\!&0.72 \, i,
\end{array}
\end{equation}
The same computations enable us to evaluate the L-function accurately,
in this case finding
$Z(5) \approx -0.03656426734$.  This $Z$-function is the bottom one graphed in Figure \ref{fig:zfunctions}.  

An important point is that the test functions need to be chosen sufficiently different from each other to make
the system well-conditioned.  This is difficult to do 
when the spectral parameters are large and is a limitation of the method.   Another
interesting point is that we are forced to work to high precision, typically to more than 50 decimal places of 
accuracy, because the system is ill-conditioned and many digits of precision are lost in the calculations.  Again for more details, see \cite{FKL}.



\subsection{Rigor and completeness} 


While our methods never actually establish the existence of L-functions, 
it is worth emphasizing that our methods are capable of producing rigorous statements.  
For example, they can prove that there are no L-functions with
spectral parameters within certain regions.  Similarly, suppose
one knew somehow from some external reason that 
a given region contains exactly one L-point.  Then 
our methods can rigorously identify small intervals guaranteed
to contains the correct coordinates of the L-point and the correct real and imaginary parts of
some initial coefficients~$a_p$. 

In practice, the step in our search process where we are most likely
to miss actual L-functions is in the initial scan of a region to find
likely candidates.    It can happen for example that search parameters
which seem to work well in one region miss L-functions in a nearby region. 
As one leaves the origin and crosses contour lines in any of 
Figures~\ref{fig:r0r0r0}, \ref{fig:r0r1r1}, and \ref{fig:Cr1plot},
we expect the percentage of L-functions that our program has
found generally decreases.  We expect that we have found almost
all L-points  up through 
the first few contours, but it is surely a much smaller percentage 
towards the end.




\section{Zero-free and L-point-free regions}
\label{sec:free}
Each transcendental L-function should be viewed as just as remarkable
an analytic object as an algebraic L-function, even though it is 
much harder to access computationally.    We begin by 
considering zero-free intervals on the critical line
of a given L-function and conclude by considering
L-point-free regions in landscapes of L-functions.
We sketch how these two topics are parallel in nature, 
as the explicit formula can be used to 
establish such regions of each type.  





\subsection{Three L-functions and their low-lying zeros} 
\label{ssec:threeL} 
\begin{figure}[htp]
\includegraphics[width=0.98\textwidth]{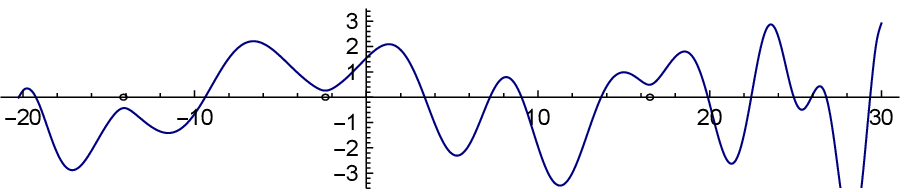}
\vskip 0.1in
\includegraphics[width=0.98\textwidth]{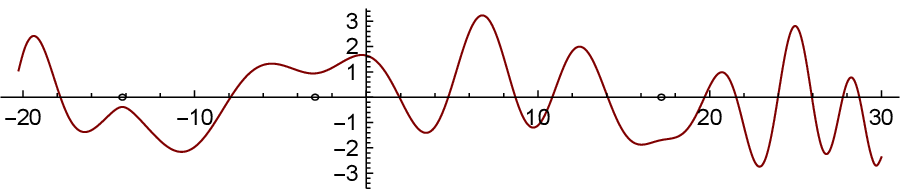}
\vskip 0.1in
\includegraphics[width=0.98\textwidth]{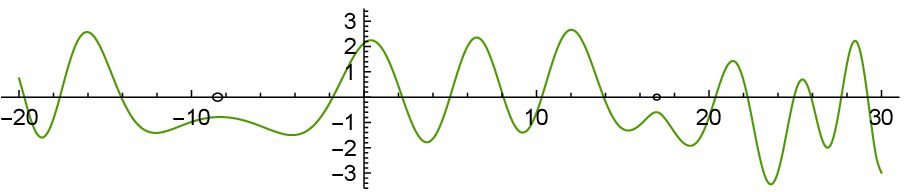}

\caption{\sf
\label{fig:zfunctions}
Three Z-functions, coming from L-functions with the following
$\Gamma$-types and L-point coordinates:  $\;\;\;\;\;\;\;\;\;\;\;\;\;\;$ \hspace{5in}
$r0r0r0$, 
$(\lambda_1, \lambda_2) \approx (14.141,2.380)$, $\;\;\;\;\;\;\;\;\;\;\;\;\;\;\;\;\;\;$ \hspace{2in}
$r0r1r1$, 
$(\lambda_1, \lambda_2) \approx (14.204,2.980)$, $\;\;\;\;\;\;\;\;\;\;\;\;\;\;\;\;\;\;$ \hspace{2in}
$r0c5,\;\;\;$ $\! (\kappa, \lambda) = (\frac52,16.976...)$. $\;\;\;\;\;\;\;\;\;\;\;\;\;\;\;\;\;\;$ \hspace{2in}
The small circles are the projections of the trivial zeros onto the critical line.}
\end{figure}

Every L-function has an associated Z-function, 
characterized up to sign by
$\abs{L(\frac12 + it)} = \abs{Z(t)}$, with
$Z(t)$ real and smooth for $t\in \R$.  
To fix the sign, we follow the LMFDB's convention
of requiring that $Z(t)$ be positive for $t$ a sufficiently
small positive number.  
Figure~\ref{fig:zfunctions} graphs the 
Z-function behind an 
L-point on each of the three figures 
from the previous section.

The Riemann Hypothesis for an L-function is equivalent to the assertion 
that the nontrivial zeros of the corresponding Z-function are real.
The plots in Figure~\ref{fig:zfunctions} make it appear that the Riemann Hypothesis is not
true for these three L-functions. Specifically, a positive local minimum 
or a negative local maximum is caused by non-real zeros of the Z-function, or,
equivalently, zeros off the critical line for the L-function.  But in fact
these local extrema do not point to counterexamples, because the Riemann Hypothesis is a statement about
the \emph{nontrivial} zeros of the L-function.  The trivial zeros occur at the poles of the $\Gamma$-factors, 
and it is these
trivial zeros which are causing the behavior which seems surprising at first glance.
This phenomenon is difficult to see for algebraic L-functions where 
all the trivial zeros occur on the negative real axis.

\begin{figure}[htp]
\includegraphics[width=0.3\textwidth]{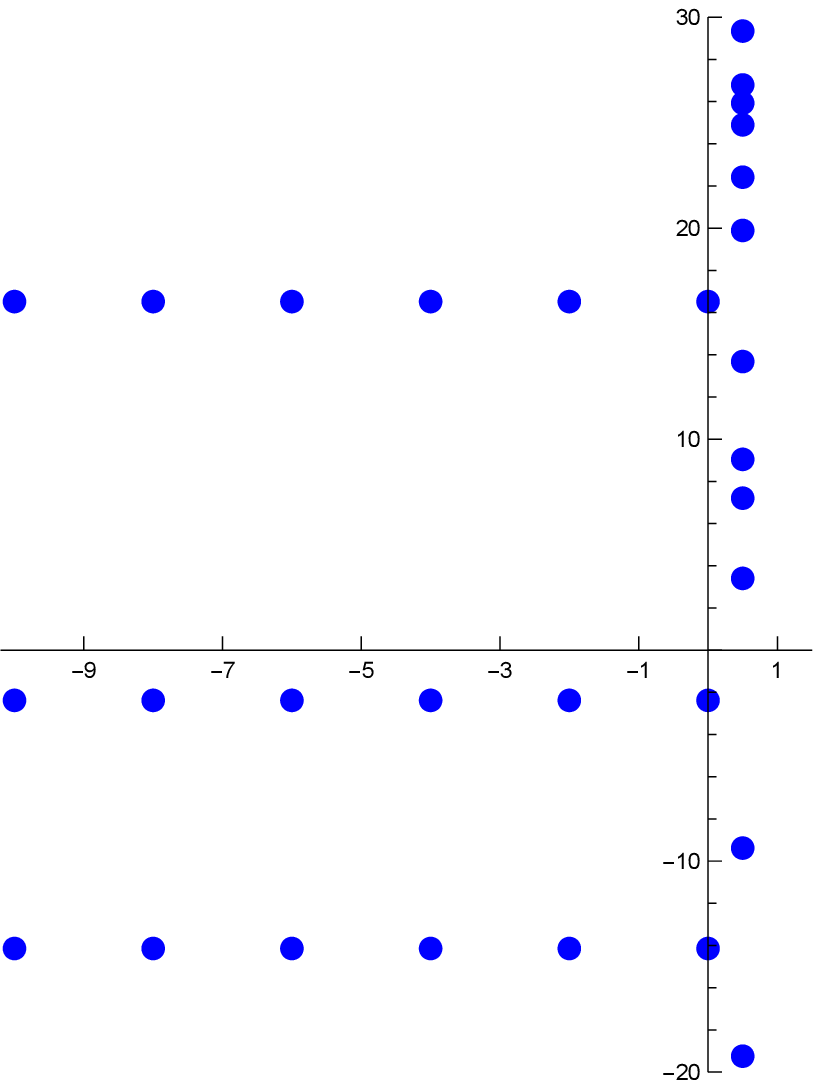}
\hfill
\includegraphics[width=0.3\textwidth]{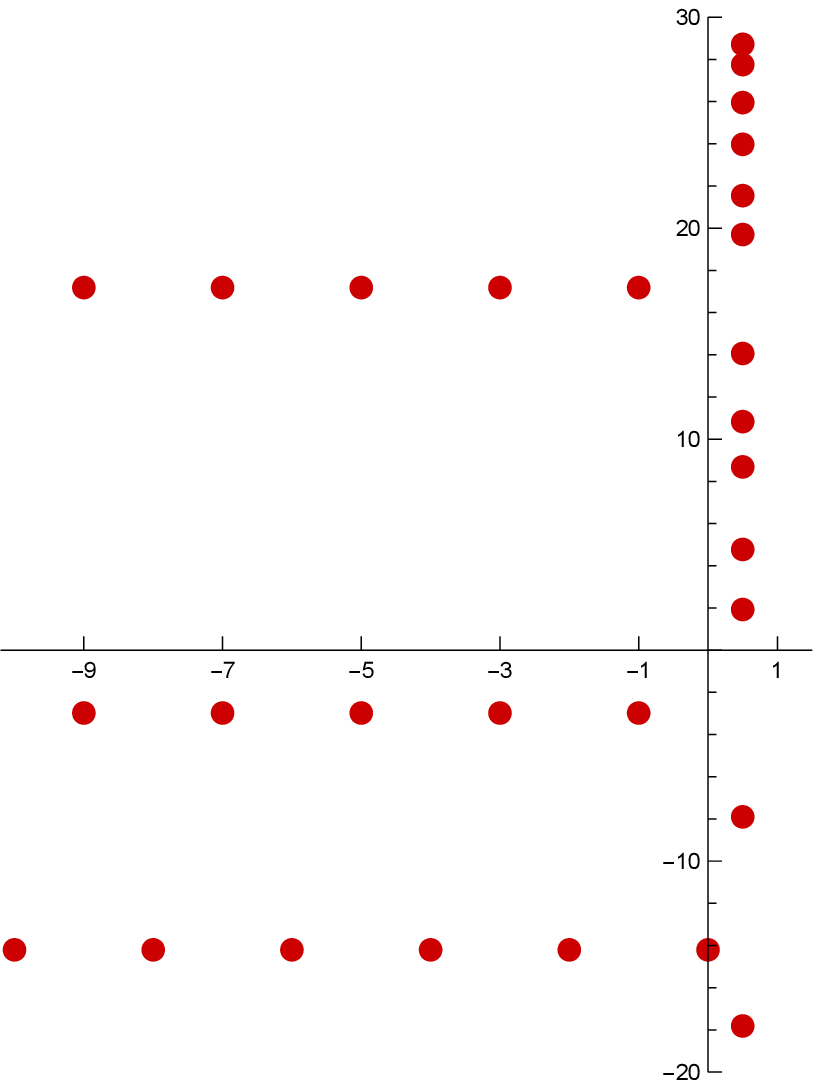}
\hfill
\includegraphics[width=0.3\textwidth]{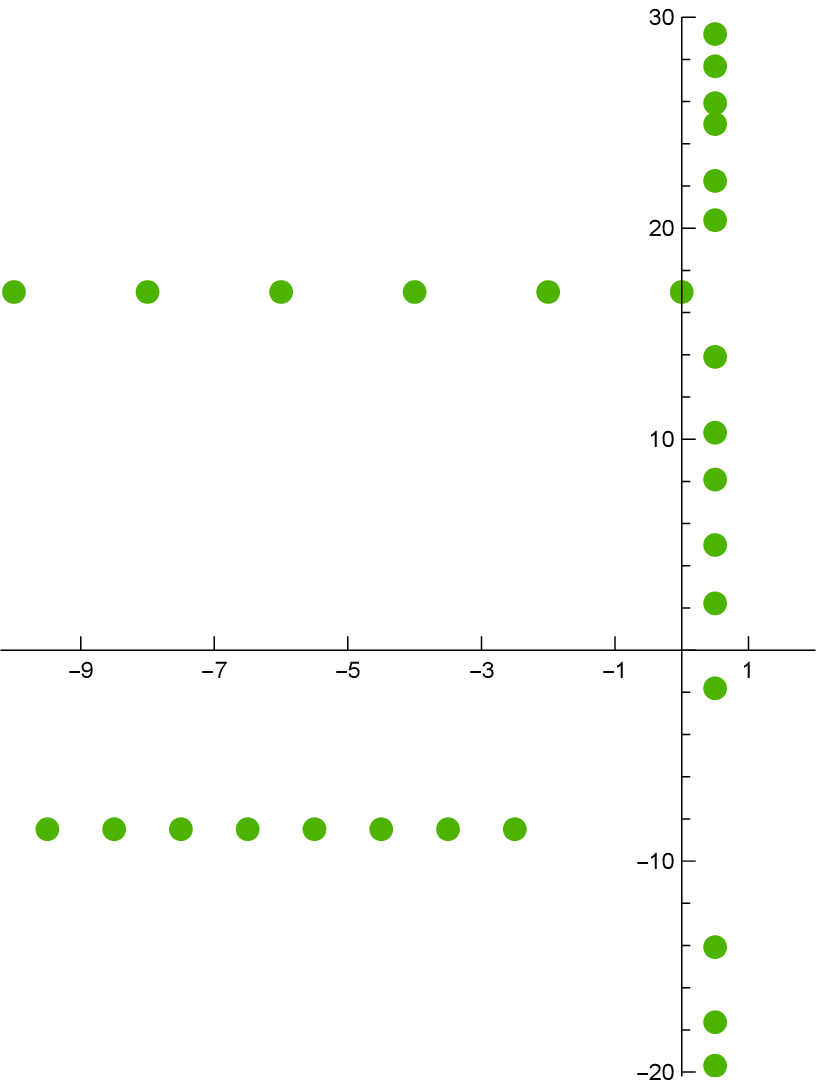}
\caption{\sf
The zeros in the complex $s$-plane of the L-functions in Figure~\ref{fig:zfunctions}.
In the region shown, the nontrivial zeros all lie on the critical line $\Re(s)=\frac12$.
} \label{fig:zeros}
\end{figure}

Figure~\ref{fig:zeros} shows how the trivial zeros
appear to ``push aside'' the nearby zeros on the
critical line.    The L-functions have been chosen to all have
trivial zeros at height approximately $17$ to facilitate 
comparison.  The trivial zeros at this height come from a gamma factor 
$\Gamma_\R(s+\delta + i \lambda)$ with $\lambda \approx -17$.  
The pushing is stronger in the 
first and third cases because the trivial zeros are one
step closer: $\delta = 0$ in these cases, while
$\delta = 1$ in the second case.   
%

\subsection{The explicit formula explains zero-free intervals}
\label{sec:explicitLfunction}  The pushing-aside phenomenon was first 
observed in~\cite{strom}.  As shown by \cite{eight} it can be 
rigorously established using Weil's explicit formula.  
This formula
%
%
in cartoon form says
\begin{equation}\label{eqn:explicit}
\sum_\gamma f(\gamma) = \hat{f}(0)\cdot\text{(log $N$)}
+ 
\int f(t)\frac{\Gamma'}{\Gamma}\text{-terms}\,dt 
+
\sum_n \hat{f}\bigl(\frac{\log n}{2\pi}\bigr)\text{coeffs} .
\end{equation}
Here $f$ is a test function, $\hat{f}$ is its Fourier transform,
$N$ is the conductor of the given L-function,  
and the sum on the left is over the zeros.  The argument in \cite{eight}
chooses $f$ to be non-negative, so the left side of 
\eqref{eqn:explicit} can be interpreted as the weighted count of the
zeros within the support of~$f$.   The test function is chosen to 
have support in an interval $I$ containing the heights of the 
trivial zeros in question.  Here the $\Gamma'/\Gamma$ 
terms on the right are negative, the effect being stronger
when the zeros are closer to the critical line.  In this way, with 
the conductor small enough as well,  one can prove 
certain intervals $I$ indeed contain no zeros.   
\mpar{Applications of the explicit formula usually come in two versions, one where the RH for all L-functions in questions is assumed, one where it isn't.  I think what we are doing here could be clarified.}


\subsection{The explicit formula and the approximate functional equation explain L-point-free regions}
\label{ssec:explicitLpoint}   A striking feature of
Figures~\ref{fig:r0r0r0}, \ref{fig:r0r1r1}, and \ref{fig:Cr1plot} 
is that there are no L-points at all in certain large regions.
Such an L-point-free region was established
for the case of Figure~\ref{fig:r0r0r0} in \cite{miller2} using the explicit  
formula \eqref{eqn:explicit}.  The L-point-free
region was enlarged in \cite{FKL} by using the
approximate functional equation \eqref{eq:formula}.  
See Figure~1 in \cite{FKL}.

The basic idea behind the proof in \cite{miller2} is that for a point $(\lambda_1,\lambda_2)$ 
in the region in question one can choose a test function $f$ taking only nonnegative values in \eqref{eqn:explicit}
so that the right side of~\eqref{eqn:explicit} is negative.  But 
the left side of \eqref{eqn:explicit}  is nonnegative and the conclusion is that no 
L-function with coordinates at $(\lambda_1,\lambda_2)$ can exist.   
This same general method has been applied in many contexts.  For 
example, simply taking $(\lambda_1,\lambda_2)=(0,0)$ one
is saying that certain L-functions coming from number
fields can only exist if conductors are large enough 
\cite{JR2}.   \mpar{I put in what I think is the most relevant reference, which happens to be a paper of mine.  
I commented out other references because I think we are over-referencing here}
The method, and the more 
powerful exclusionary techniques of \cite{FKL}, can 
be extended to establish L-point-free regions of other landscapes.  

For the discussion in \S\ref{ssec:compare2} and \S\ref{ssec:compare3} we are interested
in how the L-point-free region varies with the
parameter space.  For example, consider the 
two different parameter spaces which are superimposed in Figure~\ref{fig:Cr1plot},
for which the Plancherel measures are exactly the same.
Observationally, the lowest purple points on lines where
 $\kappa$ is integral are lower than the lowest green points
 on lines where $\kappa$ is half-integral.  Accordingly,
 we expect a larger zero-free region in $X_{+c}$ than
 in $X_{-c}$.  From the viewpoint of the explicit formula,
 the situation is similar to the contrast between the
 three cases in Figure~\ref{fig:zeros}.  Namely the factor $\Gamma_\R(s+i\lambda)$
 gives a larger negative contribution to the $\Gamma'/\Gamma$-term
 in \eqref{eqn:explicit} than the factor $\Gamma_\R(s+1+i \lambda)$ does.

\section{Coefficient space}
\label{sec:coefficient}  Here we treat  general degrees $d$ and explain how the union $X'_d$ of all
the connected parameter spaces $X'_{d_+,d_-,d_2}$ maps finite-to-$1$ to a single 
\term{coefficient space} $Y'_d = \R^{d-1}$.   
Information is lost, because of the failure of injectivity, but many parameter
landscapes, like Figures~\ref{fig:r0r0r0}, \ref{fig:r0r1r1}, and \ref{fig:Cr1plot} are combined into a single coefficient landscape,
like Figure~\ref{eucliddata}.  


\subsection{Spectral parameters as three multisets of complex numbers}
\begin{figure}[htb]
\begin{center}
\includegraphics[width=5in]{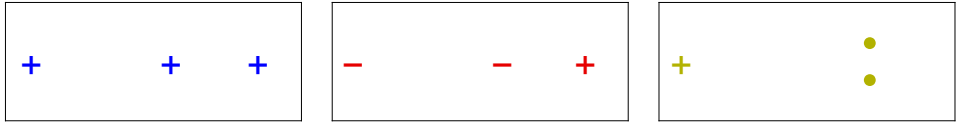}
\caption{Dot diagrams for the $\Gamma$-factors of Figure~\ref{fig:zfunctions},
drawn in the rectangle with $|\mbox{Re}(z)| \leq 20$ and $|\mbox{Im}(z)| \leq 8$.  \label{dotdiagrams}}
\end{center}
\end{figure}

An equivalent way to give spectral parameters is to give multisets $Z_+$, $Z_-$, 
and $Z_2$ of complex numbers in a slightly different way, as follows:   
\begin{align*}
\Gamma_\R(s + i \lambda)  \ \mbox{ contributes } \ & \lambda  \mbox{ to } Z_+,  \\
\Gamma_\R(s + 1+ i \lambda)  \ \mbox{ contributes } \ & \lambda  \mbox{ to } Z_-,  \mbox{ and }\\
\Gamma_\C(s+\kappa + i \lambda)  \ \mbox{ contributes } \ & \lambda + i \kappa \mbox{ and } \lambda - i \kappa  \mbox{ to } Z_2.             
\end{align*}
We call the triple $(Z_+,Z_-,Z_2)$ a {\bf dot diagram} because we think of it visually 
as illustrated by Figure~\ref{dotdiagrams}.
So the space $X_d$ of allowed spectral parameters in degree $d$ is identified with the space of dot diagrams having
$d$ points, counting multiplicities.   Repeating material from \S\ref{ssec:spaces} in the new language, $X_d$
decomposes as the disjoint union of the subspaces $X_{d_1,d_2}$ parameterizing 
dot diagrams with $d_2$ complex conjugate pairs of non-real roots.  Each $X_{d_1,d_2}$ 
then decomposes into the disjoint union of the subspaces $X_{d_+,d_-,d_2}$ parametrizing
dot diagrams with $|Z_+|=d_+$, $|Z_-|=d_-$, and $|Z_2|=2d_2$. To understand the previously-defined
topology of all these spaces one should think of the dots as being allowed to move 
horizontally in the complex plane.   Dots are allowed to cross and dots of a given type are indistinguishable. 
 The ambient spaces $X_d'$ introduced in \S\ref{ssec:spaces} correspond to letting
complex conjugate pairs of dots move arbitrarily in $\C-\R$.  

The viewpoint on the $X_{d_+,d_-,d_2}$  of the previous sections of course remains valid.  
Thus, as explained before, the closed yellow cones in Figures~\ref{fig:r0r0r0} and  \ref{fig:r0r1r1} respectively 
indicate half of $X_{3,0,0} = X_{+++}$ and $X_{1,2,0}=X_{+--}$.   Likewise one should think of
Figure~\ref{fig:Cr1plot} as a picture of the upper half of both 
$X_{+c} \subset X'_{+c}$ and $X_{-c} \subset X'_{-c}$, as each full
space $X_{\tau c}$ has the form $\frac{1}{2} \Z_{\geq 1} \times \R$.  

\subsection{Spectral parameters as factorizing real polynomials}
Our introduction of coordinates $(\lambda_1,\dots,\lambda_{d_1-1})$ 
and $(\kappa_1,\dots,\kappa_{d_2})$ on $X'_{d_1,d_2}$ in \S\ref{ssec:spaces} involved an 
ordering convention from \S\ref{ssec:terminology}.  As is often the case, it is good
to have convention-independent coordinates.  Accordingly, 
we form corresponding monic polynomials, 
\begin{equation}
f_\epsilon(x) = \prod_{z \in Z_\epsilon}(x-z) \in \R[x].
\end{equation}
The distinction between the $+$ points and the $-$ points 
plays a secondary role, and so it is moreover useful to consider 
the product,
\[
f(x) =  f_+(x)f_-(x)f_2(x) =   x^d + c_2 x^{d-2} + \cdots + c_{d-1} x +  c_d.
\]
Thus a point $(Z_+,Z_-,Z_2)$ of $X'_d$ determines
a point $(c_2,\dots,c_d)$ of $\R^{d-1}$.   The dual
of $(Z_+,Z_-,Z_2)$ is $(-Z_+,-Z_-,-Z_2)$.  In the new
coordinates, dualizing corresponds to replacing $f(x)$ by $(-1)^d f(-x)$ 
so that $c_j$ is replaced by $(-1)^j c_j$. 

Let $Y_d \subset \R^{d-1}$ be 
the image of $X_d$ and let $Y_{d_1,d_2}$ be the image of 
$X_{d_1,d_2}$.  The map 
$X_{d_1,d_2} \rightarrow Y_{d_1,d_2}$
corresponds to forgetting the distinction between the $+$ and $-$ points.
It is thus a branched cover of degree $2^{d_1}$.
Numerics are governed by a row
of Pascal's triangle, as the branched cover
 $X_{d_+,d_-,d_2} \rightarrow Y_{d_1,d_2}$
 has degree the binomial coefficient 
 $\frac{d_1!}{d_-! d_+!}$.  
 These numerics stay the same if 
 we consider the enlarged spaces 
 $X'_{d_+,d_-,d_2}$ and their images
 $Y'_{d_1,d_2}$.   
 
 \subsection{The coefficient landscape}
\begin{figure}[htb]
\begin{center}
\includegraphics[width=5in]{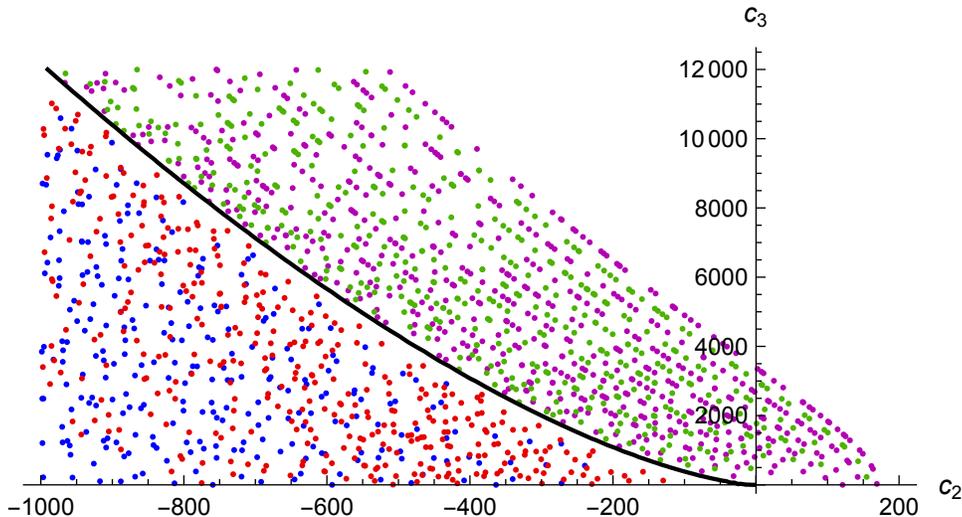}
\caption{\label{eucliddata} 
Landscape of L-functions in the coefficient plane $Y'_3 = \R^2$. }
\end{center}
\end{figure}

 Figure~\ref{eucliddata} draws a window on the upper half of the coefficient plane $Y'_3=\R^2$. 
Also included are the data points from Figures~\ref{fig:r0r0r0}, \ref{fig:r0r1r1},
\ref{fig:Cr1plot}, and their duals.   The space $Y_{3,0}$ is
the subspace where the discriminant $D=-27 c_3^2-4 c_2^3$
is positive or zero.  It is filled by blue points coming from 
the bijection $X_{+++} \rightarrow Y_{3,0}$ 
and red points from the triple cover $X_{+--} \rightarrow Y_{3,0}$.

For a real number $\kappa$,
let $C_\kappa$ be the contour curve with equation $D=-4 \kappa^6$,
so that $C_0$ is the drawn discriminant locus bounding $Y_{3,0}$.  
The space $Y_{1,1}$ is a union of curves $C_\kappa$ 
indexed by $\kappa \in \frac{1}{2} \Z_{\geq 1}$.
The maps from $X_{-c}$ and $X_{+c}$ to 
$Y_{1,1}$ are both bijective.  The curve $C_\kappa$ 
contains purple points from $X_{-c}$ or green points 
from $X_{+c}$, according to whether $\kappa$ is an integer or not. 

A key feature of Figure~\ref{eucliddata} and also Figure~\ref{2pict} is a decomposition
of $\R^{d-1}$ into two regions via the discriminant $D$, a left one where $D \geq 0$ and 
a right one where $D<0$.   In the general case, there are $\lceil d/2 \rceil$ parts, 
the images of the $X'_{d_1,d_2}$.  All the components of $Y_d$ are entirely 
inside one of these parts.  The components of $Y_d$ in the image of $X'_{d_1,d_2}$ 
all have codimension $d_2$.

\section{Plancherel measure}
\label{sec:plancherel} In the parts near the origin, where data should be at least nearly complete,
there is a noticeable difference between the density of L-points in parameter spaces,
as illustrated by Figures~\ref{fig:r0r0r0}, \ref{fig:r0r1r1}, and \ref{fig:Cr1plot}, and the density of L-points in coefficient
space, as illustrated by Figure~\ref{eucliddata}.   Here we explain
how Plancherel measure gives a clear theoretical framework which
in particular explains the visual difference.   



\subsection{The global meaning of Plancherel measure}  
\label{globalmeaning}
Plancherel measure is naturally defined via the group $\GL_d(\R)$ and plays a large role
in many purely local questions.   However, we are interested in this paper 
only in its conjectural role in counting L-functions and so will
take a more global viewpoint.   For a measurable set 
$U \subset X_d$ and a positive integer $N$, let $\cL(U,N)$ 
be the set of L-functions with $\Gamma$-factor in $U$ and 
conductor $N$.   Then the governing principle
is that one should have 
\begin{equation}
\label{planchapprox}
|\cL(U,N)| \approx \mu_d(U) \nu_d(N),
\end{equation}
for a multiplicative function $\nu_d$ on positive integers coming from
analogous Plancherel measures for the $p$-adic groups $\GL_d(\Q_p)$.  
Rigorously, one 
conjectures an asymptotic equivalence as 
the cutoff $x$ goes to infinity,
\begin{equation}
\label{planchasymp}
\sum_{N \leq x} |\cL(U,N)| \sim \sum_{N \leq x} \mu_d(U) \nu_d(N).
\end{equation}
The general asymptotic equivalence \eqref{planchasymp} 
can only hold for at most one pair of measures, so 
\eqref{planchasymp} is in particular a conjectural global characterization
of the Plancherel measure $\mu_d$.  Theorem~1 of \cite{LV}, with a special case given
in \eqref{Laplacian} below, is one 
of many strong statements in the literature supporting the 
general truth of \eqref{planchasymp}.

As a simple illustration of this formalism, let $d=1$ where
$X_{1} = X_+ \coprod X_-$.  Each summand on the right
has one element and maps bijectively to $Y_1 = \{0\} = \R^0$.   
The number $\nu_1(p^j)$ is the number of characters
$(\Z/p^j)^\times \rightarrow \C^\times$ that do not factor through the quotient group  $(\Z/p^{j-1})^\times$.
 Thus, e.g., $\nu_1(p)=p-2$ for
any prime $p$.  To give a primitive Dirichlet character 
$\chi: (\Z/N)^\times \rightarrow \C^\times$ with $N = \prod_p p^{e_p}$ is
to give an injective character 
$\chi_p : (\Z/p^{e_p})^\times \rightarrow \C^\times$ for
each~$p$.  For $N \geq 2$ the two possible values
$1$ and $-1$ arise equally often for $\chi(-1)$.   The Plancherel
 measure $\mu_1$ gives mass $1/2$ to
 both $X_+$ and $X_-$.   So in fact \eqref{planchapprox}
 holds as an equality for all $U$ and all $N \geq 2$.  
 The behavior at $N=1$ represents a doubling
 phenomenon which holds for all $d$: half the 
 $\Gamma$-factors are excluded and the others
 should asymptotically occur at twice
 their Plancherel density.

\subsection{Modified absolute values and near-constant densities}
This subsection and the next are essentially translations of presentations in the literature 
into more intuitive terms.  See, e.g., \cite{Hel} for a more Lie-theoretic presentation.  
 As a preliminary for writing down Plancherel measures on parameter spaces $X_{d_1,d_2}$ 
 in an intuitive way,
 define 
 \begin{align*}
{|t|_+} &  {= t \, \tanh \left( \frac{\pi}{2} t  \right)}, &
{|t|_-}& { = t \, \coth \left( \frac{\pi}{2} t \right)}, &
|t|_0 & = |t|.
\end{align*}
To help understand Plancherel measures on coefficient
spaces $Y_d$, restrict to $t \geq 0$ and consider the change of variables $s=t^2$, and thus $t = s^{1/2}$ and 
$dt = \frac{1}{2}s^{-1/2} ds$.  Then the differential form $t f(t) dt$ is equal to $s^{1/2} f(s^{1/2}) \frac{1}{2} s^{-1/2} ds = \frac{1}{2}f(s^{1/2}) ds$.

\begin{figure}[htb]
\begin{center}
\includegraphics[width=5in]{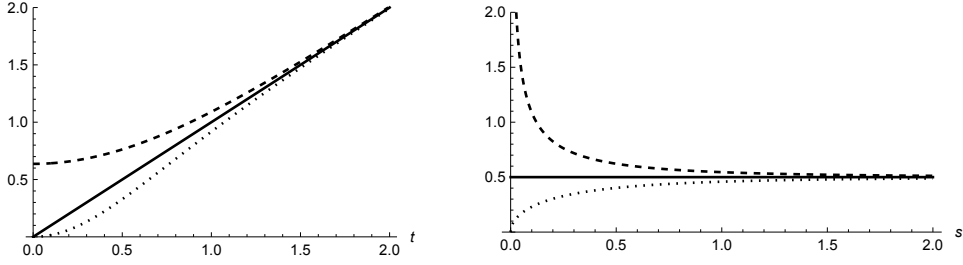}
\end{center}
\caption{\label{threeabs} Left: densities $|t|_+$ (dotted),  $|t|_-$ (dashed), and  $|t|_0=|t|$ (solid).  Right:
associated densities $\tanh(\pi \sqrt{s}/2)/2$, $\coth(\pi \sqrt{s}/2)/2$, and $1/2$.}
\end{figure}

Figure~\ref{threeabs} shows that $|t|_+$ has a double-zero at $t=0$ while $|0|_- = 2/\pi$. 
This distinction is important when comparing with data for large $N$, but in our context 
of $N=1$ it is negligible, as we will discuss further.  
More important to us is the fact, obvious from the defining formulas, that
the modifications $|t|_\epsilon$ each approach $|t|$ exponentially fast
as $|t| \rightarrow \infty$.   So the corresponding densities in the $s$-variable
likewise converge rapidly to $1/2$.

\subsection{Plancherel measure and its near-Euclidean nature}
For constants $P_d$  specified in the next subsection, the Plancherel measure is given on the 
$(d_1+d_2-1)$-dimensional part $X_{d_1,d_2}$ of $X_d$ by 
\begin{align}
\label{planchform}
\mu_d & =  \frac{P_d}{2^{d_1+d_2}}  \prod_{i<j} |z_i - z_j|_{\epsilon_i \epsilon_j} \;\; d \lambda_1 \cdots d \lambda_{d_1+d_2-1}.
\end{align}
Here the $z_i$ run over the points in the dot diagram and $\epsilon_i$ is $+$, $-$ or $0$ according
to whether $z_i$ is in $Z_+$, $Z_-$, or $Z_2$.  The factor $2^{d_1}$ in the denominator arises from the 
previously mentioned fact that the map from $X_{d_1,d_2}$ to its image $Y_{d_1,d_2}$ has degree~$2^{d_1}$.
The factor $2^{d_2}$ goes away in \eqref{planchform2} because each of the $\kappa_i$ are allowed to 
be $\frac{1}{2} \Z_{\geq 1}$, rather than just in $\Z_{\geq 1}$.  
 
   Formula \eqref{planchform} lets one supplement the vision of moving points in dot diagrams: the points tend to repel each other 
in the sense that moving any two points further apart contributes to an increase in Plancherel
density.   At close range, the distinction between the subscripts $+$, $-$, and $0$ is important, but
at long range it is not.  

In our cases in the setting $d=3$, using $\lambda_1+\lambda_2+\lambda_3=0$ in the first two cases, 
\eqref{planchform} becomes these more explicit formulas:
\begin{align*}
\mbox{On $X_{+++}$,} &&
 \mu_3 &= \frac{P_3}{8} |\lambda_1-\lambda_2|_+ \; |2 \lambda_1+ \lambda_2|_+ \; |\lambda_1+ 2 \lambda_2|_+ \; d \lambda_1 \, d\lambda_2,  \\
\mbox{On $X_{+--}$,} && \mu_3 &= \frac{P_3}{8} |\lambda_1-\lambda_2|_- \; |2 \lambda_1+ \lambda_2|_+ \; |\lambda_1+ 2 \lambda_2|_- \;  d \lambda_1 \, d\lambda_2, \\
\mbox{On $X_{-c}$ and $X_{+c}$,} && \mu_3& = \frac{P_3}{8} \kappa (4 \kappa^2 + 9 \lambda^2) \; d\lambda.
\end{align*}
The first two densities are indicated by light contour plots in the background in Figures \ref{fig:r0r0r0} and \ref{fig:r0r1r1} respectively, 
where they are near-indistinguishable.  While the third measure is supported on a union of vertical lines in 
Figure~\ref{fig:Cr1plot}, the background there nonetheless gives the contour plot of $\kappa (4 \kappa^2 + 9 \lambda^2)$
on the entire positive quadrant.  

Now consider the enlarged parameter spaces $X_d'$ obtained by allowing each $\kappa_j$ to run over 
$(0,\infty)$ as discussed before.   On each piece $X'_{d_1,d_2}$ define
 \begin{equation}
\label{planchform2}
\mu'_d =  \frac{P_d}{2^{d_1}}  \prod_{i<j} |z_i- z_j| \;\; d \lambda_1 \cdots d \lambda_{d_1+d_2-1} \; d{\kappa_1} \cdots d{\kappa_{d_2}}.
\end{equation}
A standard fact about discriminants, explained in \cite{JR} for example, then says that the push-forward of 
$\mu'_d$ to coefficient space $\R^{d-1}$ is just the Euclidean measure $P_d \, dc_2 \dots dc_d$.  
Keeping in mind that the differences between the three $|t|_\epsilon$ decay exponentially as $t$ becomes large,
 as illustrated by Figure~\ref{threeabs}, and our root distances are almost always large, 
 as illustrated by Figure~\ref{dotdiagrams}, the conclusion
is that Plancherel measure on coefficient space is extremely close to a Euclidean measure for our purposes.  

\subsection{Plancherel constants}    To get $P_d$, we compare the main theorem of \cite{LV},
as made explicit by \cite{Gold}, with a much more elementary result from \cite{JR}. 

\vspace{0.05in}
\noindent \term{Laplacian asymptotics.}
The Laplacian eigenvalue associated to an automorphic L-function with $(d_1,d_2) = (d,0)$ and spectral 
parameters $(\lambda_1,\dots,\lambda_{d})$ in the conventions of \cite{LV} is 
\begin{equation}
\label{laplacian}
\Lambda = 1 + \frac{\lambda^2_{1}+\cdots+\lambda^2_{d}}{2}.
\end{equation}
Placing an instance of $|\cL(U,1)|$ from \eqref{planchapprox} in the notation of \cite{LV}, 
let $N(T)$ be the number of such L-functions with conductor
$N=1$ and Laplacian eigenvalue $\Lambda \leq T$.  The symmetric space $M = K \setminus \PGL_d(\R)/\PGL_d(\Z)$
has dimension $m=\sum_{j=2}^d j = (d+2)(d-1)/2$.  Theorem 1 of \cite{LV}, specialized to
the group $\PGL_d(\R)$, gives the asymptotic
\begin{equation}
\label{Laplacian}
N(T) \sim \frac{\mbox{vol}(M)}{\Gamma(\frac{m}{2}+1) (4 \pi)^{m/2}} T^{m/2}.
\end{equation}
Here the volume depends on the ratio of choices of Haar measure for $\PGL_d(\R)$ and $K$.  
With the choices described in \cite[Section~3]{LV}, the volume is calculated in \cite[Theorem~1.6.1]{Gold} as 
\begin{equation}
\label{goldfeld}
\mbox{vol}(M) = \frac{d}{\pi^{m/2}} \prod_{j=2}^d \zeta(j) \Gamma\Bigl(\frac{j}{2}\Bigr).
\end{equation}

\noindent \term{A Euclidean volume.} 
On the other hand, the part $Y_{d,0}^r$ of coefficient space $\R^{d-1}$ coming from $(\lambda_1,\dots,\lambda_d)$ 
with $\lambda_1^2 + \cdots + \lambda_d^2 \leq r^2$ is just the part of $Y_{d,0}$ with $c_2 \geq -r^2/2$.   For 
example the upper half of $Y_{3,0}^{\sqrt{2000}}$ is the part of Figure~\ref{eucliddata} beneath drawn
discriminant curve $C_0$.  The Euclidean volume of $Y_{d,0}^r$, 
with respect to the measure $dc_2 \cdots dc_d$, is calculated via a Selberg integral
in \cite[Props.\ 1.1 and 1.2]{JR} to be 
\begin{equation}
\label{selberg}
\mbox{vol}(Y_{d,0}^r) = \frac{\prod_{j=2}^d \Gamma(\frac{j}{2})}{2^{d(d-1)/4} \sqrt{d} \Gamma(\frac{m}{2} + 1)} r^{m}.
\end{equation}
Note that the right sides of \eqref{Laplacian}, \eqref{goldfeld}, and \eqref{selberg} all have a slightly different form
from the right sides in the references.  We have made these small changes, using facts like the duplication
formula for the $\Gamma$-function, so as to make the simplification \eqref{cd} short.  

\vspace{.05in}
\noindent \term{Comparison.}
 If we replaced $T$ on the right side of \eqref{Laplacian} by $T-1$, the asymptotic would still be true.  
So, the desired Plancherel constant $P_d$ is the ratio of the right sides of \eqref{Laplacian} and \eqref{selberg}, 
with $T$ replaced by $r^2/2$ in conformity with \eqref{laplacian}:

\begin{eqnarray}
\nonumber P_d & = & \frac{\frac{d}{\pi^{m/2}} \prod_{\ell=2}^d \zeta(\ell) \Gamma(\frac{\ell}{2})}{\Gamma(\frac{m}{2}+1) (4 \pi)^{m/2}} (r^2/2)^{m/2} \cdot
 \frac{2^{d(d-1)/4} \sqrt{d} \Gamma(\frac{m}{2} + 1) }{ \prod_{j=2}^d \Gamma(\frac{j}{2})} r^{-m} \\ 
\nonumber  & = & \frac{\frac{d}{\pi^{m/2}} \prod_{\ell=2}^d \zeta(\ell)}{(4 \pi)^{m/2}} 2^{-m/2} \cdot
 2^{d(d-1)/4} \sqrt{d}\\ 
\label{cd} & = & \frac{{d}^{3/2}}{2^{(d+3)(d-1)/2}} \prod_{j=2}^d \frac{\zeta(j)}{\pi^j}.
\end{eqnarray}
The final formula is written so that the factor $\zeta(j)/\pi^j$ is rational for even $j$.  
Special cases of \eqref{cd} are $P_1=1$, $P_2=\frac{1}{12}$ and the evaluation \eqref{eqn:P3} of $P_3$ 
highlighted in the introduction.

\subsection{Comparison of data and theory in $d=2$.}  
\label{ssec:compare2} In this subsection
and the next, we compare three things: computed data, Plancherel
measure $\mu_d$, and the Euclidean approximation 
$\mu_d'$ to Plancherel measure.    Here we summarize the 
well-understood $d=2$ case, using it to calibrate our
expectations for our $d=3$ case.  

\begin{figure}[htb]
\begin{center}
\includegraphics[width=5in]{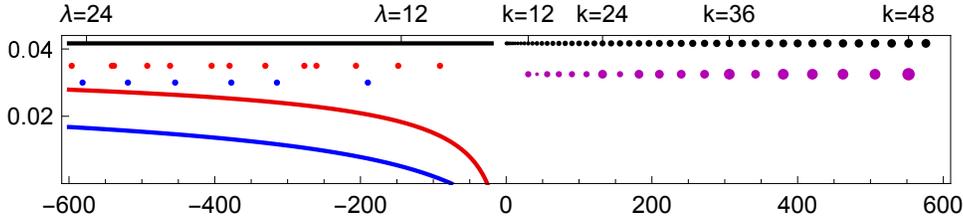}
\end{center}
\caption{\label{2pict} Plancherel measure on the $c_2$-interval $[-600,600]$
indicated in black at the top.  Colored L-points come from Maass forms on
the left and holomorphic cusp forms on the right.  }
\end{figure}

We are guided in the $d=2$ case by Figure~\ref{2pict},
which is to be compared with Figure~\ref{eucliddata} for the
$d=3$ case.   The Euclidean approximation $\mu_2' = \frac{dc_2}{12}$ 
gives mass $50$ to both the intervals
$[-600,0]$ and $[0,600]$.  But, as indicated by 
the colored L-points, taken from \cite{LMFDB},
there are only $19$ L-functions in the negative half, giving the very low percentage $38\%$.   For the positive half,
write always $k=2 \kappa+1$ and thus $\kappa = (k-1)/2$.  An L-point 
at $\kappa \in \frac{1}{2} \Z_{\geq 1}$ arises with multiplicity the dimension of the space of 
cusp forms $S_k(1)$ of weight $k$ on $\SL_2(\Z)$.  For odd $k$, this dimension
is always $0$.  For $k=12$, $14$, $16$, $18$, $20$, and $22$ the dimensions
are $1$, $0$, $1$, $1$, $1$, and $1$.  Adding $12$ to an even 
$k$ increases the multiplicity by $1$, for a general approximate
formula of $(k-1)/12$.  These multiplicities are indicated
by the areas of the L-points, with the rightmost drawn L-point coming from
$k=48$ and its multiplicity of $4$.  The total number of L-points, including multiplicity,
works out to $37$, for the still small percentage 
of $74\%$.   So we are interested here in
obtaining some understanding of why these
percentages are so small.

For $c_2>0$, it is helpful to use the function $f(k) = ((k-1)/2)^2$.  
Plancherel measure is supported on $f(\Z_{\geq 2})$.  By the case of \eqref{planchform} with no differentials it gives mass
$(k-1)/24$ to $f(k)$.  A natural cutoff is $C=f(K+1/2)$ for a positive integer 
$K$.  Then, assuming that $K$ is a multiple of $12$ for the last formula, 
\begin{align}
\label{eqn:muprime} \mu'_2([0,C]) & = \frac{f(K+1/2)}{12}  \!\!\!\!\! &&\!\!\!\!\! \!\!\!\!\!\!\!\!\! = \frac{K^2}{48} - \frac{K}{\;\,48\,\;} + \frac{1}{192}, \\
\label{eqn:mu}\mu_2([0,C]) & = \sum_{k=2}^K \frac{k-1}{24}  \!\!\!\!\!  &  &\!\!\!\!\!   \!\!\!\!\!\!\!\!\! = \frac{K^2}{48} - \frac{K}{\,\;48\;\,}, \\
\label{eqn:actual} |\cL([0,C],1)| & = \sum_{k=2}^K \dim(S_{k}(1)) \!\!\!\!\!  \!\!\!\!\!  & &\!\!\!\!\!   \!\!\!\!\!\!\!\!\!  = \frac{K^2}{48} - \frac{14 K}{\,\;48\;\,} + \frac{244}{192}.
\end{align}
So the low-percentage problem is not at all a failure of $\mu_2'$ to be a good 
approximation for $\mu_2$.  One can rewrite the discrepancy between the theories \eqref{eqn:muprime}, \eqref{eqn:mu} and the data \eqref{eqn:actual} in the $c_2$ variable.  One then gets that for $c_2>10$ the density 
\begin{equation}
\label{second00}
f_{c}(c_2) = \frac{1}{12} \left( 1 - \frac{3.25}{\sqrt{c_2}} \right)
\end{equation}
approximates the density of L-points much better than $\frac{1}{12}$ does.     
Note also that the doubling phenomenon arises here: at level $N=1$, 
odd $k$ never occur and so even $k$ occur twice as often as 
Plancherel measure would predict.  

For $c_2<0$ there is still a doubling phenomenon: at level $N=1$, dot 
diagrams of the form $+-$ cannot occur and so $++$ and
$--$ asymptotically occur twice as often as Plancherel measure would predict.  
So both $++$ and $--$ should have asymptotic average density
$1/24$. This exact constant comes from the Euclidean measure $\mu_2'$ and is graphed 
in Figure~\ref{2pict}.  The density for actual Plancherel measure $\mu$ would visually
differ from $1/24$ only in the $c_2$ interval $[-2,0]$, by the asymptotic constancy of the right half of 
Figure~\ref{threeabs}.  Five-term refinements for the counting functions with 
a very small error terms are given in  \cite[Eqs.\ 40,41]{Steil}.
Taking only three terms, the corresponding 
densities written in terms of $b=|c_2|$ are
\begin{align}
\label{secondpp} f_{++}(b) & = \frac{1}{24} \left( 1 - \frac{9 \log(b) + 3 \log(\pi^4/2)}{\pi \sqrt{b}} \right),  \\
\label{secondmm} f_{--}(b) & = \frac{1}{24} \left(1  - \frac{3 \log(b) + 3 \log(8)}{8 \pi \sqrt{b}} \right).
\end{align}
These functions are plotted in Figure~\ref{2pict}.   Their low values over the plotted region are consistent with there being 
only $6$ L-functions of type $++$ and $13$ L-functions of type $--$, rather than the $25$ of each that 
Plancherel measure alone would predict.   The slow approach of these functions
to their common asymptote is illustrated by the numerical values
 $24 f_{++}(10000) \approx 77\%$ and
$24 f_{--}(10000) \approx 89\%$.

\subsection{Comparison of data and theory in $d=3$} 
\label{ssec:compare3} We can now 
compare the coefficient landscape drawn in Figure~\ref{eucliddata} with expectations
from theory.  The Plancherel measure of any rectangle of size say $200\mbox{-by-}2000$ 
is almost exactly $(200)(2000)P_3 \approx 210$.    From the previous
subsection, we should not expect to obtain this density in the range of the figure, 
but we should expect to be closer in the region with $d_1=1$ 
to the right than we are in the region with $d_1=3$ to the left.  
 A greater density is indeed visually apparent on the right.  
 For example, the rectangle with $-300 \leq c_2 \leq -100$ and 
$2000 \leq c_3 \leq 4000$ has $119$ green and purple L-points, giving the percentage 
$119/210 \approx 57\%$.  In comparison, the rectangle with
$-500 \leq c_2 \leq -300$ and $0 \leq c_3 \leq 1000$ in
the totally real region on the left has only $82$ blue and red L-points, which is
about $39\%$ of $210$. 

We think it is possible that we have found all the L-functions in the two sample
rectangles just mentioned. If one goes closer to $(c_2,c_3) = (0,0)$, we think it is
 more likely that we have found all the L-functions.  Here the percentages are
generally lower, and we interpret that as being from secondary effects, capturing
not only that there are no L-points in the L-point-free regions
discussed in \S\ref{ssec:explicitLpoint}, but that L-points close
to these regions should be sparse.   If one goes further from $(c_2,c_3)=(0,0)$
again the experimental density tends to decrease, but now it is clear that
our searches have simply missed L-functions.  

In interpreting the difference between the blue L-points from $X_{+++}$  and the red
L-points from $X_{+--}$ on the totally real side, one needs to keep in mind that there are 
asymptotically three times as many of the latter as there are of the former.  
But even with this in mind, the data near the discriminant locus $C_0$ 
suggests that secondary effects for +++  are much stronger than
they are for $+$$-$$-$; this is expected from comparison of 
the $++$ and the $--$ cases in $d=2$.   More data would be very 
helpful in making a more detailed analysis.  Even more helpful would
be refined densities with more terms, analogous to 
\eqref{second00}, \eqref{secondpp},
and \eqref{secondmm}.


\section{Complements}
\label{sec:complements}
    We conclude by pointing to two promising future directions.  
 We describe them briefly here.  Both involve the partially 
 conjectural concept of the Sato-Tate group associated to an 
 L-function.  

\subsection{Comparison with the $p$-adic Plancherel formula}
\label{ssec:pplancherel}
    Another way to examine our collection of L-functions is to compare
 it with expectations from the spherical $p$-adic Plancherel formula.
 This formula asymptotically governs the distribution of 
 coefficients $a_p \in \mathbb{C}$ for a fixed prime $p$ and varying 
 L-functions.  There are uniform formulas for all~$d$, but we will 
 be briefer than in the last section and consider only our case~$d=3$.   
 
 \begin{figure}[htb]
\begin{center}
\includegraphics[width=5in]{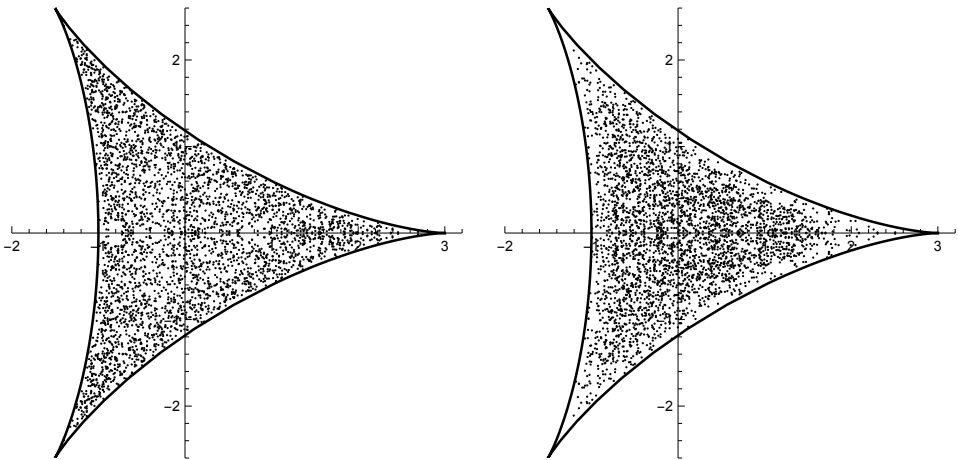}
\end{center}
\caption{\label{twovsfive} Coefficients $a_p$
from the degree~3 conductor~1  L-functions in our collection: $a_2$'s on the left and $a_5$'s on the right.} 
\end{figure} 
 
      To see some data before bringing in theory,   Figure~\ref{twovsfive} plots
the $a_2$'s and $a_5$'s of all the degree~3 conductor~1 L-functions in our collection.
Except in a few cases those L-functions are not self-dual, which explains the
reflection symmetry across the real axis.
  The Ramanujan conjecture says that each $a_p$
should agree with the trace of an element of the compact group $\SU_3$.  The conjecture holds here, as expected, because the $a_p$ indeed lie in the 
curvilinear triangle $T$ of traces.  
Writing the generic element in the complex plane as $z = x+i y$,
  the boundary of $T$ is given by the vanishing of
 \[
 f(x,y) = 27 - 18 x^2 + 8 x^3 - x^4 - 18 y^2 - 24 x y^2 - 2 x^2 y^2 - y^4.
 \]
 For a fixed non-self-dual L-function and varying $p$, the $a_p$ are conjectured
 to be equidistributed in $T$ according to the push-forward of the Haar 
  measure on $\SU_3$.  The density of this Sato-Tate measure $\mu_\infty = f_\infty(x,y) dx dy$ 
  works out to $f_\infty(x,y) = \sqrt{f(x,y)}/(2 \pi^2)$.  
 
 The $p$-adic Plancherel measures $\mu_p = f_p(x,y) dx dy$ fit together to 
 form a family indexed by $p \in [1,\infty]$, with $p=\infty$ giving the  Sato-Tate measure just discussed. 
  The density $f_1(x,y)$ is constant.  The 
 formula for general $p$ is too complicated to give here, 
 but for $p>1$ all the $f_p(x,y)$ vanish
 on the boundary, and the unique maximum is at $(0,0)$ and 
 increases with $p$.   A~precise statement which captures this phenomenon of
 increasing concentration in the middle is
 \begin{equation}
 \label{theoreticalmoment}
 \int_T f_p(x,y) (x^2+y^2) dx dy = 1 + \frac{1}{p} + \frac{1}{p^2}.
 \end{equation}
 We obtained all these statements by specializing the general definition of 
 $p$-adic Plancherel measure \cite{Macdonald} to $\SU_3$ and then carrying
 out the requisite multivariate calculus.  

Figure~\ref{twovsfive} shares qualitative features with the behavior predicted
by $p$-adic Plancherel measures.  For example, the distribution of the $a_2$'s
is hard to visually distinguish from the uniform distribution, but the distribution
of the $a_5$'s already shows marked repulsion from the boundary.  
However, quantitatively, experimental moments are far off from 
theoretical moments:
\[
{
\renewcommand{\arraystretch}{1.1}
\begin{array}{r|ccc}
p & 2 & 3 & 5 \\
\hline
\mbox{Theoretical moment from \eqref{theoreticalmoment} }& 1.75 & 1.\overline{44} & 1.24 \\
\mbox{Experimental moment from dataset} & 1.52 & 1.17 & 0.96 \\ 
\end{array}
}
\]
We again expect that this discrepancy is explained by as yet unknown secondary
terms.  

\subsection{Group-theoretical framework}  
To make our concrete points in the next two subsections cleanly, we
first bring in more formalism.  
To understand the full set $\cL$ of all automorphic L-functions it is clarifying to 
use the conjectural \term{Langlands group} $G$ and
two of its quotients:
\begin{equation}
\label{threegroups}
G \twoheadrightarrow G^{\rm alg} \twoheadrightarrow G^{\rm fin}.
\end{equation}
Here $G$ is a very large compact group containing Frobenius classes $\Fr_p$ for all primes $p$.  The set  
of irreducible representations of $G$ is identified with $\cL$ by 
requiring that a representation $\rho$ corresponding to an    
L-function $L(s) = \prod_p f_p(p^{-s})^{-1}$  satisfy $f_p(x) = \det(1-\rho(\Fr_p) x)$ for all primes $p$
not dividing the conductor.  The representations  that factor through $G^{\rm alg}$ 
correspond to L-functions in $\cL^{\rm alg}$, therefore L-functions for which the
$\lambda_j$ are all $0$.  The representations that factor through $G^{\rm fin}$ 
correspond to L-functions for which the $\kappa_k$ are also all $0$.  
 Conjecturally, $G^{\rm alg}$ has an 
independent definition as a compact version of the absolute
motivic Galois group of $\Q$.  Conjecturally, 
$G^{\rm fin} = \Gal(\overline{\Q}/\Q)$.   The
\term{Sato-Tate group} of an L-function $L$ is the
image of the corresponding representation $\rho$;
it is a compact subgroup of $\GL_d(\C)$, well-defined
up to conjugation.   In general Frobenius
elements are conjectured to be equidistributed
in the conjugacy classes of the Sato-Tate group,
as mentioned in the previous subsection
for the Sato-Tate group $\SU_3$.

If the Sato-Tate group of an L-function has $c$ components, then 
the component group conjecturally comes from a Galois
field of degree $c$.   Any prime dividing the discriminant of this field
divides the conductor of the L-function.   So requiring
$N=1$ forces Sato-Tate groups to be connected.  
Also the abelianizations of the three groups \eqref{threegroups} 
are expected to coincide with the abelianization $\hat{\Z}^\times$
of $\Gal(\overline{\Q}/\Q)$.  So, at $N=1$, tori are also
ruled out as possible Sato-Tate groups. 

 For $d=3$ and $N=1$ there are only two possibilities.  
 The Sato-Tate group is $\SO_3$ in the special case that
 the L-function is self-dual and $\SU_3$ otherwise.  
 This simple dichotomy allowed us to ignore the group-theoretical
 framework in previous sections.  But when one allows
 general $N$, there are infinitely many possible $G$, 
 with the finite groups arising in \S\ref{ssec:algebraicL} illustrating some 
 possibilities.  To regain some simplicity in discussing
 the landscapes for $(d,N)$, one can focus on 
 \term{generic} L-functions.  By definition, these 
 are L-functions whose Sato-Tate groups 
 contain $\SU_d$, with $d$ being the degree.   
 One expects that for any fixed $d$, $100\%$ of
 L-functions are generic asymptotically, for say
 a ball centered at the origin in coefficient space $\R^{d-1}$ 
 of increasing radius.  

\subsection{Generic algebraic L-functions are hard to find for $d \geq 3$}
\label{ssec:algebraicpart}
We now support a main theme of the introduction, that algebraic L-functions
are relatively rare for degrees $d \geq 3$.  The classical degrees
do not point in the right direction at all: in degree $1$, all 
L-functions are algebraic. For $d=2$, an L-function is algebraic if and only
if $c_2 \geq 0$ as illustrated by Figure~\ref{2pict}; so Plancherel measure
says in a strong sense that half of the L-functions are algebraic and 
half are not.



  
The set $Y_3^{\rm alg}$ is just the set of 
points $(c_2,c_3) = (\kappa^2,0)$ for 
$\kappa \in \Z_{\geq 0}$.  We impose genericity for 
the rest of this subsection, thereby excluding the parameter $\kappa=0$.  
Computations with the cohomology of $\SL_3(\Z)$ in \cite{AP} say
that there are no generic L-functions in $\cL^{\rm alg}_{3,1}$ with $N=1$ and 
$\kappa < 120$.  To see some algebraic L-functions, we increase the
conductor, like we did in \S\ref{ssec:algebraicL}.  At $\kappa=1$, cohomological computations
with congruence subgroups in \cite{GKTV}
show the existence of dual pairs of generic L-functions
for the conductors $89$, $106$, $116$, $128$, $160$, $205$, $212$, and $221$.
An L-function with L-point $(\kappa^2,0)$ 
for $\kappa \in \Z_{\geq 1}$ should 
come from a motive with nonvanishing 
Hodge numbers $h^{2 \kappa,0}=h^{\kappa,\kappa}=h^{0,2 \kappa}=1$.  
For $\kappa=1$, a one-parameter family of 
such motives was found in the cohomology of surfaces in \cite{GT}.  
A~dual pair of members of this family numerically matches
the dual pair of automorphic L-functions with conductor
$128$.   For $\kappa>1$ there is a qualitative
difference, as Griffiths transversality prevents 
motives from moving in families,
and so algebraic geometry can only produce them 
one dual pair at a time.

In summary,  generic algebraic L-functions 
in degree $d=3$ are known only for $\kappa=1$.
For $d \geq 4$, the only examples we know of come from families of threefolds in \cite{scholten,DV}, with
nonvanishing Hodge
numbers $h^{3,0}=h^{2,1}=h^{1,2}=h^{0,3}=1$; the smallest conductor seems to be $2^9 3^9$,
so matching automorphic calculations are out of reach.

%

\subsection{Comparison in the self-dual setting} 
For a more subtle comparison, consider now only L-functions in $\cL$ with real 
coefficients.  According to the group-theoretic
framework, they come in two types, orthogonal
and symplectic, the latter arising only for 
$d$ odd.  In the spirit of the Frobenius-Schur indicator,
we use $\tau \in \{+,-\}$ to distinguish them,
with $+$ for orthogonal and $-$ for symplectic.  

The theoretical context set up in this paper has direct analogs in 
both of these two parallel settings.  Thus there are sets
of L-functions $\cL^\tau_{d,N}$ mapping
to disconnected parameter spaces $X^\tau_{d}$ 
of possible $\Gamma$-factors
which in turn map to coefficient 
spaces $Y^\tau_{d} \subset \R^r$. Here $r = \lfloor d/2 \rfloor$ arises as the
dimension of a maximal torus in the
compact group $\Sp_d$ or $\Orth_d$.  
There is again a Plancherel measure $\mu_\tau$ 
on $Y_d^\tau$ and again a canonical Euclidean approximation $\mu'_\tau$ 
on $\R^{r}$. Again $\mu_\tau$ should govern the distribution of
L-points.   A~difference is that $Y^\tau_d$ contains 
components of all the possible dimensions 
$0$, \dots, $r$.   So Figure~\ref{2pict} for $\SU_2 = \Sp_2$ is
a good guide.  Like in the $d=2$ case, the isolated points are
algebraic and each has positive Plancherel measure.   
Another difference is that algebraic geometry provides
many large families of source motives for the
algebraic points of many of the larger-dimensional components.  For
example, the $g(g+1)/2$-dimensional family
of dimension $g$ abelian varieties maps 
to the algebraic point of a $\lfloor g/2 \rfloor$-dimensional 
component of $Y_{2g}^-$.   

As mentioned in the introduction, some transcendental
L-functions representing the symplectic case with
$(d,N) = (4,1)$ are already in the LMFDB.  
The $d=4$ cases with somewhat larger conductors
should be within computational reach.




\end{document}